\documentclass{elsart}

\usepackage{epsfig}
\usepackage{latexsym}
\usepackage[all]{xy}
\usepackage{amssymb}
\usepackage{amsmath}
\usepackage{amsxtra}

\DeclareMathOperator{\Gal}{Gal}

\DeclareMathOperator{\Det}{Det}

\DeclareMathOperator{\Tr}{Tr}
\DeclareMathOperator{\Norm}{N}

\newcommand{\bfQ}{\mathbb Q}

\newcommand{\bfP}{\mathcal P}

\newcommand{\p}{\mathfrak p}
\newcommand{\q}{\mathfrak q}
\newcommand{\bB}{\mathfrak B}
\newcommand{\bA}{\mathfrak A}
\newcommand{\F}{\mathcal F}

\newcommand{\R}{\mathbb R}
\newcommand{\Ba}{\mathbb B}
\newcommand{\Ta}{\mathbb T}

\newcommand{\Q}{\mathbb Q}
\newcommand{\C}{\mathbb C}
\newcommand{\Z}{\mathbb Z}

\newcommand{\De}{\mathfrak D}
\newcommand{\D}{\mathcal {D}}
\newcommand{\M}{\mathcal {M}}
\newcommand{\Siegel}{\mathfrak h}
\newcommand{\lat}{\mathfrak L}
\newcommand{\disc}{\mathfrak d}

\newcommand{\A}{\mathcal{A}}

\newcommand{\B}{\mathcal{B}}
\newcommand{\scrO}{\mathcal {O}}

\journal{Journal of Number Theory}
\begin{document}
\begin{frontmatter}
\title{A formula for the central value of certain Hecke L-functions}
\author{Ariel Pacetti}
\ead{apacetti@mate.dm.uba.ar}
\address{Departamento de Matem\'atica, Pabell\'on I - Ciudad
  Universitaria. (1428)\\ Buenos Aires, Argentina}
\thanks{This paper is the author's Ph.D. thesis presented at the
  University of Texas on August 2003. He would like to thank
  Dr. Fernando Rodriguez-Villegas and the University of Texas at Austin.} 
\thanks{The author would like to thank and dedicate this paper to
  Joana Terra for all her support.} 

\begin{abstract}              
  Let $N \equiv 1 \bmod 4$ be the negative of a prime,
  $K=\Q(\sqrt{N})$ and $\scrO_K$ its ring of integers. Let $\D$ be a
  prime ideal in $\scrO_K$ of prime norm congruent to $3$ modulo $4$.
  Under these assumptions, there exists Hecke characters $\psi_{\D}$
  of $K$ with conductor $(\D)$ and infinite type $(1,0)$. Their L-series
  $L(\psi_\D,s)$ are associated to a CM elliptic curve $\bA(N,\D)$
  defined over the Hilbert class field of $K$. We will prove a
  Waldspurger-type formula for $L(\psi_\D,s)$ of the form
  $L(\psi_\D,1) = \Omega \sum_{[\A],I} r(\D,[\A],I) m_{[\A],I}([\D])$ where the sum is
  over class ideal representatives $I$ of a maximal order in the
  quaternion algebra ramified at $|N|$ and infinity and $[\A]$ are
  class group representatives of $K$. An application of
  this formula for the case $N=-7$ will allow us to prove the non-vanishing
  of a family of L-series of level $7|D|$ over $K$.

\end{abstract}
\end{frontmatter}

\section{Introduction}
Given an elliptic curve $E$ over $\Q$, and a fundamental discriminant
$D$, a formula of Waldspurger relates the value of $L(E
\otimes D,1)$, the twist of $E$ by $D$, with the coefficients of a
$3/2$ modular form (see \cite{Waldspurger}). The purpose of this work
is to get a formula for quadratic twists of a family of elliptic
curves with complex multiplication not defined over the rationals.

Given an imaginary quadratic field $K$ the theory of complex
multiplication (see \cite{Shimura2}) gives a relation between elliptic
curves with CM given by an order of $K$ and L-series associated to
Hecke characters $\psi$ on $K$. The simplest case is when
$K=\Q(\sqrt{N})$ with $N \equiv 1 \bmod 4$ the negative of a prime and
$\psi$ is a character of conductor $\sqrt{N}$. In this case the
L-series corresponds to a CM elliptic curve $\bA(N)$ studied by Gross
in \cite{Gross}, defined over $H$, the Hilbert class field of $K$. A
formula for the central value of $L(\psi,1)$ was given by Villegas in
\cite{Villegas}.

In this paper we will study the central value of the L-series
corresponding to the CM elliptic curves $\bA(N,\D)$, given by twists
of $\bA(N)$ by the quadratic character of conductor $\sqrt{N}\D$ where
$\D$ is a prime ideal of $K$ prime to $\sqrt{N}$ and with prime norm
congruent to $3$ modulo $4$. If we denote $h$ the class
number of $K$, the prime ideal $\D$ has $h$ Hecke
characters $\psi_\D$ of conductor $\D$ associated to it. The relation
between the L-series of $\bA(N,\D)$ and 
$L(\psi_\D,s)$ is given explicitly by :
$$L(\bA(N,\D)/H,s) = \prod_{\psi_\D} L(\psi_\D,s) L(\overline{\psi_\D},s)$$
where $H$ is the Hilbert class field of $K$ and the product is over
the $h$ Hecke characters associated to $\D$ (see \cite{Gross} formula
(8.4.4) and Theorem 18.1.7). If we define $\bB$ be the Weil restriction
of scalars of $\bA(N,\D)$ to $K$, then $\bB$ is a CM abelian variety, and
$L(\bA(N,\D)/H,s) = L(\bB/K,s)$.

Let $B$ be the quaternion algebra ramified at $|N|$ and
infinity. Given an element $x \in B$ we denote $\Norm(x):= x \bar x$
its norm and $\Tr(x) := x + \bar x$ its trace. To the ideal $\D$ and
an element $[\A]$ of $Cl(\scrO_K)$ we will associate a maximal order
$O_{[\A],[\D]}$ in $B$ depending only on $[\A]$ and the class of $\D$. If $\{I\}$ are
representatives for left $O_{[\A],[\D]}$-ideals, the main theorem (Theorem
\ref{main-teo}) gives the formula $L(\psi_\D,1) = \Omega \sum_{[\A],I}
r(\D,[\A],I) m_{[\A],I}([\D])$ where the sum is over the ideals
$\{I\}$ and ideal representatives of $\scrO_K$, $\Omega$ is a period,
$r(\D,[\A],I)$ is a rational integer and the numbers
$m_{[\A],I}([\D])$ are algebraic integers.

The paper consists of four chapters besides the introduction. In the
second chapter we give the basic definitions and derive a 
first formula for the value of the L-series at $1$ (following Hecke's
work on L-series, see \cite{Hecke}). Later we relate Theta functions
of quadratic forms to Theta functions on the Siegel space.
In the third chapter we introduce the period $\Omega$ and using
Shimura's theory in Complex Multiplication we compute the field where
the algebraic integers $m_{[\A],I}$ belong to. 
In the fourth chapter we study the problem of deciding whether two
points in the Siegel space are equivalent or not in our specific
case. For this purpose we introduce quaternion algebras, and relate
special points with left $O_{[\A],[\D]}$-ideals.
In the last chapter we study in detail the case when the
class number of $K$ is one. In this case the elliptic curve $\bA(N)$
is defined over $\Q$ and the numbers $m_{[\A],I}$ turn out to be
rational integers. In the case $N=-7$ using the fact that the
quaternion algebra has class number $1$ for maximal ideals, we prove
that the CM elliptic curves $\bA(N,\D)$ defined over $K$ have a
non-vanishing L-series for all primes $\D$. \\
We finish this work with a remarkable relation between the numbers
$m_{[\A],I}$ and the coordinates of the eigenvector of the modular form
associated to $\bA(N)$ represented in the Brandt matrices of level
$N^2$.

The author would like to thank the referee for suggesting a better
organization of this work.

\section{L-series}

\subsection{L-series definition}

Given a number field $K$, we will denote $\scrO_K$ its ring of
integers, $Cl(\scrO_K)$ its class group and $h$ its class number.\\
Let $N \equiv 1 \bmod 4$ be the negative of a prime, $N \neq -3$ and $K :=
\Q(\sqrt{N})$. Let $D \equiv 1 \bmod 4$ be
the negative of a prime such that the ideal generated by $D$ splits
completely in $K$, i.e. $(D) = (\D) \bar{(\D)}$. We will denote
$L:=\Q(\sqrt{D})$. Since the rings $\scrO_K \slash \D$ and $\Z \slash
|D| \Z$ are isomorphic we define $\varepsilon_\D$ by: 
\[
\xymatrix{ (\scrO_K \slash \D)^\times \ar[rr]^{\varepsilon_\D} \ar[rd] & & \pm 1\\ &
(\Z \slash |D| \Z)^\times \ar[ru]_{\left( \frac{}{|D|} \right)}\\ } \]
where $\left ( \frac{}{|D|} \right )$ is the Kronecker symbol. The
character $\varepsilon_\D$ induces a Hecke character $\psi_{\D}$ on
principal ideals by $\psi_{\D}(\langle \alpha \rangle) =
\varepsilon_{\D}(\alpha) \alpha$.

\begin{prop} The character $\psi_\D$ on principal ideals is well defined.
\end{prop}

\noindent {\bf Proof.} Since $1$ and $-1$ are the only units in $K$,
we must check that $\varepsilon_{\D} (\alpha) \alpha = -
\varepsilon_{\D} (-\alpha) \alpha$.  This follows from the fact that
$\varepsilon_{\D}$ is multiplicative and $|D| \equiv 3 \bmod 4$, hence
$\varepsilon_{\D} (-1) = -1$ $\square$\\ 
Given $\D$ let $\sigma_\D$ denote an element in $Gal(K^{\text{ab}}/K)$
corresponding to $\D$ via the Artin-Frobenius map, where
$K^{\text{ab}}$ denotes the abelian closure of $K$. We can define
$\varepsilon_\D$ in a different way:
\begin{prop} If $\alpha \not \in \D$, $\varepsilon_\D(\alpha) =
  (\sqrt{\alpha})^{\sigma_\D-1}$
\label{Artin-Tate}
\end{prop}
\noindent {\bf Proof.} It is clear that $\sqrt{\alpha}^{\sigma_\D} =
\xi_2 \sqrt{\alpha}$ where $\xi_2 = \pm1$. By definition given
$\De$ an ideal of $\bar K$ lying above $\D$, $\sigma_\D$ satisfies
$\xi_2 \sqrt{\alpha} = \sqrt{\alpha}^{\sigma_D} \equiv
\sqrt{\alpha}^{|D|} \bmod \De$. But $\sqrt{\alpha}^{|D|} = \alpha
^{\frac{|D|-1}{2}} \sqrt{\alpha}$ hence $\alpha^{\frac{|D|-1}{2}}
\equiv \xi_2 \bmod \D$. In particular $\varepsilon_\D(\alpha) =
\varepsilon_\D(\xi_2) = \xi_2$ since $|D| \equiv 3 \bmod 4$. $\square$

The character actually depends of the choice of $\D$ (i.e. we have one
character associated to $\D$ and another one associated to
$\bar{\D}$). Abusing notation $\psi$ will denote the
character associated to $\D$ if it makes no confusion.\\
The character $\psi$ defined on principal ideals extends to $h$ Hecke
characters on $I(\scrO_K)$ the set of ideals of $\scrO_K$. We fix an
extension once and for all and we call it $\psi$. Then $\psi :
I(\scrO_K) \longrightarrow T_{\psi}$, where $T_{\psi}$ is a 
degree $h$ field extension of $K$.

\begin{defn}
The L-series associated to $\psi$ is 
\begin{equation}
L(\psi,s) := \sum_{\A} \frac {\psi( \A)}{{N\A}^s}
\end{equation}
 where the sum is over all ideals $\A$ of $\scrO_K$.
\end{defn}

By Hecke's work we know that $L(\psi,s)$ extends to an analytic
function in the upper half plane, and satisfies the functional
equation:
$$
\left(\frac{2 \pi}{\sqrt{ND}}\right)^{-s} \Gamma(s) L(\psi,s) = w_{\psi} \left(\frac{2
\pi}{\sqrt{ND}}\right)^{s-2} \Gamma(2-s) L(\bar{\psi},2-s)
$$
\noindent where $w_{\psi}$ is the root number. The character $\psi$
defines a weight $2$ 
modular form given for $z$ in the upper half plane by $f_{\psi}(z) =
\sum_{\A} \psi({\A}) e^{2 \pi i z N\A}$, which has level $ND$. The
root number is given by
$w_{\psi}=f_{\psi}(\frac{i}{\sqrt{ND}}) /
\overline{f_{\psi}(\frac{i}{\sqrt{ND}})}$.\\
\begin{prop} Let $\alpha$ be a generator of $\D^h$. 
The root number in the functional equation for $\psi_{\D}$ is $w_{\psi} = \xi_2
  \left(\frac{2}{|N|}\right) i \frac{\alpha}{|\alpha|}$, 
  where $\xi_2$ is $-1$ if $2$ is ramified in $K(\sqrt{\alpha
  \sqrt{N}})$ and $1$ if not.
\end{prop}

\noindent{\bf Proof.} See \cite {Buhler-Gross} proposition 10.6, page
20. This is equivalent to saying that if $\alpha$ is the generator of
$\D^h$ such that $K(\sqrt{\alpha \sqrt{N}})$ is a quadratic extension
of $K$ of conductor $\sqrt{N} \D$ then $w_{\psi} = -
\left(\frac{2}{|N|}\right) i \frac{\alpha}{|\alpha|}$.  $\square$

The characters $\psi$ are associated to a CM elliptic curve $\bA(N,\D)$
defined over $H$, the Hilbert class field of $K$, by the formula:
$$L(\bA(N,\D)/H,s) = \prod_{\psi_\D} L(\psi_\D,s) L(\overline{\psi_\D},s)$$
\noindent See \cite{Gross} formula ($8.4.4$) and Theorem $18.1.7$.

\subsection{Choosing characters in a consistent way}

Let $\D$ and $\D'$ be prime ideals of $K$ as before (i.e. they have
prime norm congruent to $3$ modulo $4$). While extending the Hecke
character $\psi_{\D}$ to $I(\scrO_K)$ we get a field extension
$T_{\psi_\D}$. If we extend the Hecke character associated to $\D'$ in
an arbitrary way, the image of both characters will lie in different
fields. There is a natural way of defining a Hecke character
$\psi_{\D'}$ associated to $\D'$ such that $\psi_{\D'}(I(\scrO_K))
\subset T_{\psi_\D}$. Any ideal of $K$ raised to the $h$-power is
principal, hence for all ideals $\A$ prime to $\D \D'$ we define:
\begin{equation}
\psi_{\D'}(\A) = \psi_{\D}(\A)
\frac{\varepsilon_{\D'}(\A^h)}{\varepsilon_{\D}(\A^h)}
\label{consistency}
\end{equation}
There is some abuse of notation on this definition since although $\A^h$ is
principal, it has two generators $\alpha$ and $-\alpha$. 
But $\varepsilon_{\D}(-\alpha) = - \varepsilon_{\D}(\alpha)$ and
$\varepsilon_{\D'}(-\alpha) = - \varepsilon_{\D'}(\alpha)$ hence the
quotient is well defined.

\begin{prop} There exists a Hecke character associated to $\D'$ taking
  values in $T_\psi$ and defined as above on ideals prime to $\D \D'$.
\end{prop}
\noindent{\bf Proof.} We start by proving that the character defined
above is a Hecke character on ideals prime to $\D \D'$. If $\A$ is
principal, say $\A = \langle \alpha \rangle$, then $\psi_{\D'}(\alpha)
= \varepsilon_{\D}(\alpha) \alpha
\frac{\varepsilon_{\D'}(\alpha)^h}{\varepsilon_{\D}(\alpha)^h}$. Since
$h$ is odd, and $\varepsilon$ takes the values $\pm 1$, we get that
$\psi_{\D'}(\alpha) = \varepsilon_{\D'}(\alpha) \alpha$, hence it is a
Hecke character.

Let $\q$ be a prime ideal in the same equivalence class as $\D$ and
prime to $\D \D'$ (there exists such an ideal by the Tchebotarev
density theorem), say $\q \beta = \D$. Then $\psi_{\D'}(\D) =
\psi_{\D'}(\q \beta) = \psi_{\D'} (\q) \psi_{\D'}(\beta) =
\psi_{\D'}(\q) \varepsilon_{\D'}(\beta) \beta$. In this way we can
extend the character to all ideals prime to $\D'$ and clearly this is
well defined, taking values in $T_{\psi}$. $\square$\\
%% \\ Given a prime ideal $\p$ prime to $\D$, we will denote $\psi_{\p}$ the Hecke
%% character associated to $\p$. 
From now on given two different characters $\psi_\D$ and $\psi_{\D'}$
we will always assume that they are chosen in a consistent way.

Given a quadratic imaginary field $\Q[\sqrt{-d}]$ we denote $w_d$ the
number of units in its ring of integers. For $z \in \Siegel$, we recall the definition:
$$\eta(z) = e^{2 \pi i z/24} \prod_{n=1}^{\infty} (1 - e^{2 \pi i n
z})$$ While choosing ideal class representatives $\{[\A]\}$ for $K$ we
will assume they are prime to the ideal $(6)$ and that they are
written as $\A = \langle a, \frac{b+\sqrt{N}}{2} \rangle$ with $b
\equiv 3 \bmod 48$. We define $\eta(\A) :=
\eta(\frac{b+\sqrt{N}}{2a})$.  Our main theorem is the following:

\begin{thm} Given $\D$ a prime ideal of $K$ of prime norm congruent to
  $3$ modulo $4$ let $\psi_\D$ be a Hecke character as before. Let $B$
  be the quaternion algebra over $\Q$ ramified at $|N|$ and
  infinity. For each ideal class representative $[\A]$ of 
  $K$ there exists $O_{[\A],[\D]}$ a maximal order in $B$ such that:
\begin{equation}
L(\psi_\D,1)=\frac{2 \pi}{w_{|D|} \sqrt{|D|}} \eta(\bar \D)\eta(\scrO_K)
\left(\sum_{[\A]} \sum_{I} r(\D,[\A],I)m_{[\A],I}([\D])\right)
\label{11}
\end{equation}
\label{main-teo}
where $\{I\}$ is a set of left $O_{[A],[\D]}$-ideal representatives,
$r(\D,[\A],I) \in \Z$ and $m_{[\A],I}([\D])$ are algebraic integers
lying in a finite field extension of $\Q$ (see Diagram 1).
\end{thm}
The term $\Omega = \frac{2 \pi}{w_{|D|} \sqrt{|D|}} \eta(\bar
 \D)\eta(\scrO_K)$ on $(\ref{11})$ corresponds to a period of the
 abelian variety $\bB$ and the number $r(\D,[\A],I)$ is counting some
 special points 
%on each class $([\A],I)$ 
with a $\pm1$ weight (see Section 4.3 for details).  The
 rest of this paper will be a constructive proof of Theorem
 \ref{main-teo}.
        
\subsection{Computing the L-series value at $1$}

%% On this section we will derive a formula for
%% $L(\psi,1)$. 

Given $\A$ an ideal of $K$ , we will denote $[\A]$ its class in the
class group. We can decompose the L-series  as
\begin{equation}
L(\psi,s) = \sum
_{[\A]} \sum_{\B \sim \A} \frac {\psi( \B)}{{N\B}^s}
\end{equation}
\begin{prop} All integral ideals equivalent to $\A$ are of the form $c
\A$ for some $c \in \A^{-1}$.
\end{prop}
\noindent{\bf Proof.} easy to check. $\square$\\
Since the only units in $\scrO_K$ are $1$ and $-1$,
$$\sum_{\B \sim \A} \frac {\psi(
\B)}{{N\B}^s} = \frac{1}{2} \sum_{c \in
\bar{\A}} \frac{ \psi(c) \psi(\A)}{ \psi(N \A)} \frac{ {N
\A}^s}{ {Nc}^s} = \frac{1}{2} {N \A}^s \frac{\psi(\A)}{\psi(N
\A)} \sum _{c \in \bar{\A}} \frac{ \psi(c)}{ {Nc}^s}$$
Since $\psi$ is multiplicative $\psi(\A)
\psi(\bar{\A}) = \psi(N\A)$, then $\frac
{\psi(\A)}{\psi(N\A)} = \frac{1} {\psi(\bar{\A})}$. 
Using the fact that $N \A = N \bar{\A}$ it follows that $ \sum_{\B
\sim \A} \frac{\psi(\B)}{{N \B}^s} = \frac{1}{2} \frac{ {N
\bar{ \A}}^s} {\psi(\bar{\A})} \sum_{c \in \bar{\A}}
\frac{\psi(c)}{ {Nc}^s}$ and we can write the L-series as:
\begin{equation}
L(\psi,s) = 
\frac {1}{2} \sum_{[\A] \in Cl(\scrO_K)} \frac{ 
{N\A}^s}{ \psi(\A)} \sum _{c \in \A} \frac{c \varepsilon_{\D}(c)}{{Nc}^s}
\end{equation}
Without loss of generality, we may assume that $\A = a \Z + \frac{
b+\sqrt{N}}{2} \Z$ and $\D = |D| \Z + \frac{ b+\sqrt{N}}{2} \Z$, hence
$\A \D = a|D| \Z + \frac{ b+\sqrt{N}}{2} \Z$ (see \cite{Villegas}  \S
2.3 page 552). If $c \in
\A$ then $c = ma + n \frac{b+ \sqrt N}{2}$, and 
$\varepsilon_{\D}(c) = \varepsilon_{\D}(ma + n \frac{b+ \sqrt
N}{2})$. Since $n \frac{b+ \sqrt N}{2} \in \D$, 
$\varepsilon_{\D}(c) = \varepsilon_{\D}(a) \varepsilon_{\D}(m) =
\varepsilon_{\D}(N \A) \varepsilon_{\D}(m)$. We will denote $z_{\A}$ 
the point $ \frac{b + \sqrt{N}}{2a}$ (respectively $z_{\D}$ the point
$\frac{b+\sqrt{N}}{2|D|}$ and $z_{\A \D}$ the point
$\frac{b+\sqrt{N}}{2a|D|}$ ). Also we denote by $\sum'$ the sum
removing the zero element (or zero vector depending on the context). We have:
\begin{equation}
L(\psi,s) = \frac {1}{2} \sum_{[\A] \in Cl(\scrO_K)} \frac{ 
{N\A}^{1-s} \varepsilon_{\D}(N\A)}{ \psi(\A)} {\sum _{m,n \in \Z}} ' \, 
\frac{\varepsilon_{\D}(m) (m + z_{\A \D} |D| n)}{\Norm(m+z_{\A \D} |D|n)^s}
\label{main-for}
\end{equation}

If we change $m$ by $-m$ in the sum, since $\varepsilon_{\D}(-1) = -1$,
the term in the inner sum can be written as
$\frac{\varepsilon_{\D}(m)}{(m+(-\bar{z}_{\A \D}) |D|n){\left| m+
(-\bar{z}_{\A \D}) |D |n \right|}^{2s-2}}$, where the point
$-\bar{z}_{\A \D}$ is in the upper half plane.  This sum is related to
Eisenstein series that we define below:

\begin{defn} Let $p$ be a prime integer and $\varepsilon(m) := \left(
\frac{m}{p} \right)$. We define the Eisenstein series associated to
$\varepsilon$ by $E_1(z,s) = \sum _{m,n \in \Z} '
\frac{\varepsilon(m)} {(m+zpn){\left|m+zpn\right|}^{2s}}$. 
\end{defn}

\noindent By (\ref{main-for}) taking $p = |D|$ we get the relation:

\begin{equation}
L(\psi,s) = \frac{1}{2}  \sum_{[\A] \in Cl(\scrO_K)} \frac{ 
{N\A}^{1-s} \varepsilon_{\D}(N\A)}{ \psi(\A)} E_1(-\bar z_{\A \D},s-1)
\label{einsestein-decomposition}
\end{equation}
$E_1(z,s)$ turns out to be a modular form of weight $1$ with a
character. We need to compute its value at $s=0$ for a point $z$ in
the upper half plane. 
This was done by Hecke and its value (given in formula (\ref{l-serie})) can be found on
\cite{Hecke} (formulas $(26)$ and $(27)$ p.475).  For the reader
convenience we re-derive the formula. \\
The series of $E_1(z,s)$ converges only for $\Re(s) >
\frac{3}{2}$, but it can be analytically continued to the whole plane
and satisfy a functional equation. We will compute its value at $s=0$
using Hecke's trick. Since $\varepsilon$ is a character of conductor
$p$, we break the sum over $m$ as:
\begin{equation}
E_1(z,s) = {\sum_{m\in \Z}}' \frac{\varepsilon(m)}{m} + 2\sum _{n =1} ^{\infty} \sum_{r \bmod p}
\varepsilon(r) \sum _{m \in \Z}  \frac{1}{( zpn +r +mp) {|zpn +r +
mp|}^{2s}} 
\end{equation}

\noindent and dividing the last sum by $p^{2s+1}$ we get:
\begin{equation}
E_1(z,s) = 2L(\varepsilon,s)+ 2\sum _{n =1 } ^{\infty} 
\sum_{r \bmod p}  \frac{\varepsilon(r)}{p^{2s+1}}  \sum _{m \in \Z}
\frac{1}{( \frac{zpn +r}{p} +m) {\left| \frac{zpn +r}{p} +m \right|}^{2s}} 
\label{eins}
\end{equation}

\noindent For $z$ in the upper half plane we define: $$H(z,s) = \sum_{m
\in \Z} \frac{1}{(z+m) {|z+m|}^{2s}}$$

\begin{lem}
Let $z = x + iy$ be a point in the upper half plane, then:
\label{shim}
$$
\sum_{m=-\infty} ^{\infty} (z+m)^{-(s+1)}(\bar{z}+x)^{-s} =
\sum_{n=-\infty}^{\infty} \tau_n(y,s+1,s) e^{2 \pi i n x}
$$

where $\tau_n(y,s+1,s)$ is given by:

$\tau_n(y,s+1,s) \frac{i \Gamma(s+1) \Gamma(s)}{{(2\pi)}^{2s+1}}= \left \{ \begin{array}{cl}
n^{2s} e^{-2\pi n y} \sigma(4 \pi n y, s+1, s) & (n>0)\\
|n|^{2s} e^{-2\pi |n| y} \sigma(4 \pi |n| y, s, s+1) & (n<0)\\
\Gamma(2s) {(4 \pi y)}^{-2s}  & n=0
\end{array} \right.$

and $\sigma(y,\alpha,\beta) = \int_0 ^{\infty} (t+1)^{\alpha-1}
t^{\beta-1} e^{-yt} dt$
\end{lem}
\noindent{\bf Proof.} This is Lemma $1$ page $84$ \cite {Shimura} $\square$\\
The right side of lemma \ref{shim} equality converges for any $s>0$, so we can
compute the limit when $s$ tends to $0$ of $\tau_n(y,s+1,s)$ in the
different cases:

\smallskip

\noindent $\bullet$ {\bf Case $n=0$:} $\lim_{s \rightarrow 0} \frac{{(2 \pi)}^{2s+1}}{i
\Gamma(s+1)} \frac{ \Gamma(2s)}{ \Gamma(s)} {(4 \pi y)}^{-2s} = -i \pi$

\smallskip

\noindent $\bullet$ {\bf Case $n<0$:} $\lim_{s
\rightarrow0} \frac{ {(2 \pi)}^{2s+1}}{i \Gamma(s+1) \Gamma(s)}
{|n|}^{2s} e^{2 \pi |n| y} \int_0 ^{\infty} {(t+1)}^{s-1} t^s e^{-4
\pi |n| y t} dt =0$

\smallskip
\noindent $\bullet$ {\bf Case $n>0$:} $\lim_{s
\rightarrow 0} \frac{{(2\pi)}^{2s+1}  n^{2s}}{i \Gamma(s+1)} e^{-2 \pi
n y}  \frac{1}{\Gamma(s)}  \int_0^{\infty} {(t+1)}^s t^{s-1} e^{-4 \pi n
y t} dt$.

\noindent We just need to compute $\lim _{s \rightarrow 0}
\frac{1}{\Gamma(s)} \int_0 ^1 {(t+1)}^s t^{s-1} e^{-4 \pi n y t}
dt$. Doing integration by parts:
$$ \int_0 ^1 {(t+1)}^s t^{s-1} e^{-4 \pi n y
t}dt = \frac{2^s e^{-4 \pi n y}}{s} - \int_0^1 t^s (t+1)^{s-1}
e^{-4\pi nyt} dt -$$ $$ -\frac{1}{s}\int_0^1 t^s(t+1)^s e^{-4\Pi nyt}(-4 \pi
nyt) dt$$
The function $\Gamma(z)$ has a simple pole at $z=0$ with
residue 1. Dividing the integral by $\Gamma(s)$ and taking the limit
when $s$ tends to zero we get: 
\begin{equation}
\lim_{s \rightarrow 0} \tau_n(y,s + 1,s) = -2 \pi i e^{- 2 \pi n y}
\end{equation}

\noindent We just prove:

\begin{lem}
$\lim_{s \rightarrow 0} H(s,z) = -\pi i - 2 \pi i \sum_{n=1}^{\infty}
q^n$
\label{q-exp}
\end{lem} 
\noindent Equation (\ref{eins}) can be written as
$$E_1(z,s) = 2 L(\varepsilon,s) + 2 \sum_{n =1} ^{\infty} \sum_{r \bmod p}
\frac{\varepsilon(r)}{p^{2s+1}} H(\frac{zpn+r}{p},s)$$
which by lemma \ref{shim} is the same as:
$$E_1(z,s) = 2 L(\varepsilon,s) + 2 \sum_{n=1} ^{\infty} \sum_{r \bmod p}
\frac{\varepsilon(r)}{p^{2s+1}} \sum_{k \in \Z} \tau_k(yn,s+1,s)
e^{2 \pi i k (\frac{xpn+r}{p})}$$
Let $G(\varepsilon):= \sum _{r\bmod p} \varepsilon(r) \xi_p^{r}$ be
the Gauss sum associated to the quadratic character $\varepsilon$.
Let $\xi _p =e^{\frac{2 \pi i}{p}}$. If we take the limit as $s$ tends
to zero and use lemma (\ref{q-exp}) in the inner sum we get:
$$\sum_{r \bmod p} \frac{\varepsilon(r)}{p} (- \pi i -2 \pi i 
\sum_{k=1}^{\infty} q^{nk} \xi_p^{rk}) = - \frac{2 \pi i}{p}
G(\varepsilon) \sum_{k=1} ^{\infty} \varepsilon(k) q^{nk}$$
If $p$ is congruent to $3$ modulo $4$ it is a well known result that 
$G(\varepsilon) = i \sqrt{p}$, then:
\begin{equation}
\lim_{s \rightarrow 0} E_1(z,s) = 2 L(\varepsilon,1) + \frac{4
\pi}{\sqrt{p}} \sum_{n=1}^{\infty} \left(\sum_{d | n} \varepsilon(d)
\right) q^n 
\label{l-serie}
\end{equation}
Applying this to equation (\ref{einsestein-decomposition}) (with $p =
|D|$)  we get the value of $L(\psi,1)$. \\
We will write this number in terms of theta functions so as to 
relate the value for different ideals $\D$. 
Let $\B$ be any ideal of $L$. For $z$ in the upper half plane, we
define $\Theta_{\B} (z) = \sum_{\lambda \in \B} e^{2 \pi i z \frac{N
\lambda}{N \B}} = 1 + \sum_{n=1}^{\infty} r_{\B}(n) q^n$ where
$r_{\B}(n)$ is the number of elements $\lambda \in \B$ of norm $n
N\B$. Clearly if two ideals of $L$ are equivalent, their theta
functions are the same.

\begin{lem}
Let $w_{|D|}$ be the number of roots of unity in $L$, and $z$ a point in the
upper half plane. Then $\frac{w_{|D|}\sqrt{|D|}}{4\pi} E_1(z,0) = \sum_{[\B] \in
Cl(\scrO_L)} \Theta_{[\B]}(z)$
\label{theta-decomp}
\end{lem}

\noindent{\bf Proof.} We need to check that the $q$-expansion on both
sides is the same. The constant term on the right side is $h$,
the class number of $\bfQ(\sqrt{D})$. On the left side we have
$\frac{L(\varepsilon,1) w_{|D|}\sqrt{D}}{2\pi}$ which by the class number
formula is $h$. Since the constant term is the same, we can apply the
Mellin transform on both sides. Dividing by $w_{|D|}$ we need to prove the equality:
\begin{equation}
\sum_{n=1}^{\infty} \frac{\sum_{d | n} \varepsilon(d)}{n^s} =
\frac{1}{w}\sum_{[\B] \in Cl(\scrO_L)} \sum_{n=1}^{\infty} \frac{r_{[\B]}(n)}{n^s}
\label{clnf}
\end{equation}
\noindent Given a number field $L$ the zeta function associated to it is: 
$$\zeta _L (s) = \sum_{\B} \frac{1}{N \B ^s}$$ 
where the sum is over all integral ideals of $L$. It follows easily from the
definition that $\zeta_L(s) = \frac{1}{w} \sum_{[\B] \in Cl(\scrO_L)}
\sum_{n=1}^{\infty} \frac{r_{[\B]}(n)}{n^s}$ which is the right hand side of
(\ref{clnf}).

It is a classical result that $\zeta_L(s) = \zeta(s)
L(\varepsilon,s)$ (see for example \cite {Washington} Theorem 4.3,
page 33), then
$\zeta_L(s) = \left(\sum_{n=1}^{\infty} \frac{1}{n^s} \right) \left(
\sum_{m=1}^{\infty} \frac{\varepsilon(m)}{m^s} \right)$ which is the left
hand side of (\ref{clnf}) $\square$\\
Note that $-\bar{z}_{\A \D} = z_{\bar{\A} \bar{D}}$, hence by
equation (\ref{einsestein-decomposition}) and lemma \ref{theta-decomp}
we get:
$$
L(\psi,1) =\frac{2\pi}{w\sqrt{|D|}} \sum_{[\A] \in Cl(\scrO_K)}
\frac{\varepsilon_\D(N\A) }{ \psi(\A)} \sum_{[\B] \in Cl(\scrO_L)}
\Theta_{\B}(z_{\bar{\A} \bar{\D}})$$
By equation (\ref{consistency})
$\psi_{\bar \D}(\A) = \psi_\D(\A) \varepsilon_{\bar \D}(\A^h)
\varepsilon_{\D}(\A^h)= \psi_\D(\A) \left(\frac{N\A}{|D|}
\right)^h$. Since $h$ is odd it follows that
$\frac{\varepsilon_\D(N\A)}{\psi_\D(\A)} = \frac{1}{\psi_{\bar
\D}(\A)}$.
\begin{thm} The value at $s=1$ of $L(\psi,s)$ is given by:
$$
L(\psi,1) =\frac{2\pi}{w_{|D|}\sqrt{|D|}} \sum_{[\A] \in Cl(\scrO_K)}  \sum_{[\B] \in Cl(\scrO_L)}
\frac{\Theta_{\B}(z_{\A \bar{\D}})}{\psi_{\bar \D}(\bar \A)}
\label{main-for2}
$$
\end{thm}

\subsection{Theta functions in several variables}

The goal now is to write the identity of theorem \ref{main-for2} in
terms of theta functions in two variables so as to relate the
$L$-function values for different primes $D$. Given an element $(\vec{z} ,
\Omega)$ in
$\C^2$x$\Siegel_2$ (the Siegel space of dimension $2$), the
generalized theta function is defined by
$$\theta(\vec{z} , \Omega) = \sum_{\vec{n} \in \Z^2} \exp(\pi i
\vec{n}^t \Omega \vec{n} + 2 \pi i \vec{n}^t . \vec{z})$$ 
It satisfies
a functional equation for the group $\Gamma_{12}$ (following Igusa
notation), which is defined to be: $ \alpha = \left(
\begin{array}{cc}
A&B\\
C&D\\
\end{array} \right) \in Sp_{2g}(\Z)$ such that $A^tC$
and $B^tD$ have even diagonal. In particular:
\begin{equation}
\theta(\vec{0} , - (Q \tau)^{-1})  = 
\sqrt{det(Q)} \, (-i) \tau \theta(\vec{0}, Q \tau)
\label{third-case}
\end{equation}
\begin{equation}
\theta(\vec{z},
\Omega +B) = \xi_{\alpha} \theta(\vec{z}, \Omega)
\end{equation} where $Q$ and $B$ are symmetric, integral and even
diagonal two by two matrices, $Q$ corresponds to a positive definite
quadratic form, $\tau$ is a point in the upper half plane and $\xi_{\alpha}$ is a root of
unity.  (see \cite{Mumford}, \S 5, page 189).

$L$ is an imaginary quadratic field, so given an ideal $\B$ of
$Cl(\scrO_L)$ we can associate to it a quadratic form of discriminant
$D$ via the group isomorphism between $Cl(\scrO_L)$ and \{equivalence
classes of quadratic forms of discriminant $D$\}.
More specifically, given a quadratic form of discriminant $D$, say
$[a,b,c]$ where $b^2-4ac = D$, we associate the ideal $\langle a,
\frac{b + \sqrt{D}}{2} \rangle$; conversely given any primitive
ideal (i.e. not divisible by any rational integer greater
than $1$) $\B$, we can chose a pair of generators of the form
$\B = \langle a, \frac{b + \sqrt{D}}{2} \rangle$, and associate to
it the quadratic form $[a, b, c]$ where $c = (b^2 - D)/(4a)$. We will
denote $Q_{\B}$ the matrix
$\left( \begin{array}{cc}
2a&b\\
b&2c\\
\end{array} \right)$ associated to the quadratic form $[a,b,c]$.

Let $\B$ be a primitive ideal representing a class in $Cl(\scrO_L)$,
 say $\B = \langle a, \frac{b
+\sqrt{D}}{2} \rangle$ with $a = \Norm(\B)$. If $\alpha \in \B$ then it
can be written uniquely as $\alpha = ma + n \left(\frac{b
+\sqrt{D}}{2} \right)$. Hence $\Norm(\alpha) = a(am^2 + mnb + n^2
\frac{b^2 -D}{4a})$ and
\begin{equation}
\Theta_{\B}(z) = \sum_{(m,n) \in \Z^2} \exp{\left[ \pi i z (m,n)
\left( \begin{array}{cc}
2a&b\\
b&2c\\
\end{array} \right)
\left( \begin{array}{cc}
m\\
n\\
\end{array} \right)\right]}
\end{equation}
Since $z \in \Siegel$ and $Q_{\B}$ is symmetric, $zQ_{\B} \in
\Siegel_2$. Hence $\Theta_{\B}(z) = \theta(\vec{0}, zQ_{\B})$. So we
can rewrite the main formula of theorem \ref{main-for2} as:
\begin{equation}
L(\psi,1) =\frac{2\pi}{w_{|D|}\sqrt{|D|}} \sum_{[\A] \in Cl(\scrO_K)} 
 \sum_{[\B] \in Cl(\scrO_L)}
\frac{\theta(\vec{0},z_{\A \bar{\D}} Q_{\B})}{ \psi_{\bar \D}(\bar{\A})}
\label{main-for3}
\end{equation}

\section{Normalization of the Theta function}

Given a point $z_{\A \D}$, we define the normalizer:
$$\Upsilon(z_{\A \D}) := 
\eta(\D) \eta(\scrO_K) \psi_{\D}(\bar{\A})$$
Then the main formula (\ref{main-for3}) can be written as:
\begin{equation}
L(\psi_\D,1) =\frac{2\pi}{w\sqrt{|D|}} \left( \sum_{[\A] \in Cl(\scrO_K)} 
 \sum_{[\B] \in Cl(\scrO_L)}
\frac{\theta(\vec{0},z_{\A \bar{\D}} Q_{\B})}{ \Upsilon(z_{\A
 \bar\D})} \right) \eta(\bar\D) \eta(\scrO_K)
\label{main-for4}
\end{equation}
We are interested in studying the number:
$n_{\A,\B,\bar\D}=\theta(\vec{0}, z_{\A \bar{\D}} Q_{\B}) /
\Upsilon(z_{\A \bar{\D}})$. The normalizer $\Upsilon$ is chosen so as
to make $n_{\A,\B,\bar\D}$ an algebraic integer as we will see later.

\subsection{Complex Multiplication}
Let $\F_M$ be the field of all modular functions of level $M$ whose
$q$-expansion at every cusp has coefficients in $\Q(\xi_M)$ where
$\xi_M$ is any primitive $M$-th root of unity. Let
$K(M)$ denote the ray class field of $K \bmod M$, and for a prime
ideal $\p$ in $K$ relatively prime to $M$ (say of norm $p$),
$\sigma(\p)$ denotes the Frobenius automorphism of $K(M)/K$
corresponding to $\p$.

Following Stark's notation if $A$ is an integral matrix of determinant
relatively prime to $M$, we denote $f \circ A$ the action of $A$ on
$f$ which is characterized by the two
properties:
\begin{itemize}
\item $(f \circ A)(z) = f(Az)$ if $A \in Sl_2(\Z)$
\item $(f \circ A)(z) = \sigma_d(f) (z)$ if $A =
  \left(\begin{array}{cc} 1&0\\ 0&d\end{array} \right)$ where
  $\sigma_d \in \Gal(\Q(\xi_M)/\Q)$ is defined by $\sigma_d(\xi_M) =
  \xi_M^d$. We extend this action to $f$ by acting on the coefficients
  of the q-expansion at infinity.
\end{itemize}

\begin{thm} \label{cm} let $f(z)$ be in $\F_M$ and suppose that $(p) = \p
\bar{\p}$ in $K$ where $p$ is a rational prime such that
$(p,NM)=1$. Suppose that $\A = [\mu, \nu]$ is a fractional ideal of
$K$ with $\vartheta = \mu / \nu$ in $\Siegel$ and let $B
\binom{\mu}{\nu}$ be a basis for $\bar{\p}\A$. Then $f(\vartheta)$
is in $K(M)$ and $f(\vartheta)^{\sigma(\p)} = [f \circ
(pB^{-1})](B \vartheta)$.

If in addition $f$ is analytic in the interior of $\Siegel$ and has
algebraic integer coefficients in its $q$-expansion at every cusp,
then $f(\vartheta)$ is an algebraic integer.
\end{thm}

\noindent{\bf Proof.} This is Theorem 3 of \cite {Stark} page 213. $\square$

\begin{prop} Following the previous notation,
$\theta(\vec{0},\frac{z}{a|D|} Q_{\B}) / \eta(\frac{z}{|D|}) \eta(z)$
is in $\F_{24aD^2}$.
\label{q-exp2}
\end{prop}

For the proof we need the elementary result:

\begin{lem} if $f(z)$ is a modular form of weight $k$ and level $N$
and $D$ is a positive integer then $f(\frac{z}{D})$ is a modular
form of weight $k$ and level at most $ND$.
\end{lem}

\medskip

\noindent {\bf Proof of proposition \ref{q-exp2}.} Let $\B$ be the
ideal $\B := \Z a + \Z \frac{b+\sqrt{D}}{2}$. Then the quadratic form
associated to $\B$ is $[a,b,c]$ with $b^2-4ac=D$ and the matrix of the
bilinear form is $\left( \begin{array}{cc}
2a&b\\
b & 2c \end{array} \right)$ . The theta series
$\theta_\B$ is the theta series associated to this quadratic form hence
it has level $|D|$, weight $1$ and a character
$\epsilon(d)=\left(\frac{D}{d}\right)$ (see \cite{Ogg} Theorem 20, page VI-$25$). 
Using the previous lemma, we have that $\theta_\B(\frac{z}{a|D|})$
is a modular form of weight $1$ and level $aD^2$. 

The eta function is a modular form of weight $1/2$ and level $24$ ,
then $\eta(\frac{z}{|D|})$ has weight $1/2$ and level $24|D|$ , so the
product of the two eta functions has weight $1$ and level
$24|D|$. Hence the quotient has weight $0$ and level at most
$24aD^2$. We do not need a sharp estimate of the $q$-expansion, hence
the minimum level is not important.

From the $q$-expansion of the functions $\theta_\B$, and $\eta$ it is
clear that the $q$-expansion at infinity of
$\theta(\vec{0},\frac{z}{a|D|} Q_{\B}) / \eta(\frac{z}{|D|}) \eta(z)$
is in $\Q(\xi_{24aD^2})$, hence we just need to check this condition
at the other cusps. For that purpose we will study the $q$-expansion
of each form separately.

Since the theta function $\theta_\B$ is a modular form for
$\Gamma_0(|D|)$, there are just two inequivalent cusps which may be
taken to be $0$ and $\infty$. One transformation that send infinity to
zero is given by the matrix $S = \left( \begin{array}{cc}
0&1\\
-1&0 \end{array} \right)$ sending $z$ to $-1/z$.

The functional equation (\ref{third-case}) reads:
\begin{equation}
\theta\left(\vec 0,Q_\B^{-1} (-1/z)\right) = det(Q_\B )^{1/2}(-i)z
\theta(\vec 0,Q_\B z)
= \sqrt{|D|}(-i)z \theta(0,Q_\B z)
\end{equation} 

Since $Q_\B^{-1} = \text{Adj}\,(Q_\B) /|D|$, replacing $z$ by
 $z/|D|$ we get
\begin{equation}
\theta\left(\vec 0,\text{Adj}\,(Q_\B) (-1/z)\right) = (-i)z / \sqrt{|D|} \,
\theta(\vec 0,Q_\B
z/|D|)
\label{third}
\end{equation}

Replacing $Q_\B$ by its adjoint matrix, we see that the $q$-expansion
at $0$ includes a $4$-th root of unity and the square root of $|D|$ (the
$z$ factor actually cancels out a factor coming from the eta
function). Since $\sqrt{D} \in \Q(\xi_{D})$, the $q$-expansion of
$\theta(0,Q_\B)$ has coefficients in $\Q(\xi_{8D})$ at all cusps. 
Replacing $z$ by $z/a|D|$ we add at most ($aD^2$)-th roots of unity to the
$q$-expansions, hence the $q$-expansion of $\theta(0,\frac{z}{a|D|}
Q_\B)$ has coefficients in $\Q(\xi_{24aD^2})$ at all cusps.

We will use the following explicit version of the transformation
formula for $\eta$, which can be found in \cite {Villegas} page 560:

\begin{lem} \label{ville} Let 
$\left( \begin{array}{cc}
\alpha& \beta\\
\gamma & \delta \end{array} \right) \in Sl_2(\Z)$ with $\gamma$ even,
$\delta$ positive (and odd), and $\tau \in \Siegel$. Then

\begin{equation} \eta \left( \frac{\alpha \tau + \beta}{\gamma \tau +
\delta} \right) = \binom{\gamma}{\delta} e_{24}(\kappa) \sqrt{\gamma
\tau + \delta} \eta(\tau)
\label{eta-eq}
\end{equation}

where $\kappa = 3(\delta-1) + \delta(\beta-\gamma)-(\delta^2-1) \gamma
\alpha$.
\end{lem}

For any matrix in $\Gamma_0(2)$, the modular form
$\eta$ changes by a $24$-th root of unity, hence
its $q$-expansion at the equivalent cusps modulo $\Gamma_0(2)$ have
coefficients in $\Q(\xi_{24})$ and the $q$-expansion of
$\eta(\frac{z}{|D|})$ has coefficients in $\Q(\xi_{24aD^2})$. But
modulo $\Gamma_0(2)$ there are just two inequivalent cusps which may
be taken to be zero and infinity also.
The eta function satisfies the functional equation $\eta(-1/z) =
\sqrt{z/i} \, \eta(z)$. Hence its $q$-expansion at zero has
coefficients in $\Q(\xi_{8})$ and $\eta(\frac{z}{|D|})$ certainly has a
$q$-expansion with coefficients in $\Q(\xi_{24aD^2})$ at zero.
$\square$

\subsection{Field of definition}

\begin{thm} The number 
$\theta(\vec{0},z_{\A \bar \D} Q_{\B}) / \eta(z_{\bar \D})
\eta(\scrO_K)$ is an algebraic integer in $H$, the Hilbert class field of $K$. 
\label{def}
\end{thm}
\noindent {\bf Proof.} The eta function does not vanish in the
upper half plane so we can apply Theorem \ref{cm} and $\theta(\vec{0},\frac{z_0}{a|D|} Q_{\B})
/ \eta(\frac{z_0}{|D|}) \eta(z_0)$ is an algebraic integer in $F$ (some field
extension of $K$ containing $H$) where $z_0 = \frac{b+\sqrt{N}}{2}$
corresponds to the ideal $\scrO_K$.

Let $g(z) := \theta(\vec{0},\frac{z}{a|D|} Q_{\B}) /
\eta(\frac{z}{|D|}) \eta(z)$. Given an
element $\sigma$ of $\Gal(F/K)$ by complex multiplication theory there
exists a prime ideal $\p$ in $K$ such that $\sigma = \sigma_\p$, where
$\sigma_\p$ is the element in $\Gal(F/K)$ corresponding to $\p$ via
the Artin-Frobenius map. We want to prove that the quotient is in $H$ hence we take
$\p$ to be principal and using the Tchebotarev density theorem we may
assume that $\p \bar \p$ is prime to $\A$ ,$\bar \D$ and the ideal $(6)$. 

Since $\bar \p$, $\A$ and $\bar \D$ are prime to each other, it easily seen
that $b$ can be chosen such that $\bar{\p} = \langle \frac{b+\sqrt{N}}{2} ,p \rangle$,
$\A = \langle \frac{b+\sqrt{N}}{2} , a \rangle$ ,$\bar \D = \langle
\frac{b+\sqrt{N}}{2}, |D| \rangle$ and $\scrO_K = \langle
\frac{b+\sqrt{N}}{2} , 1\rangle$. Let $z_0$ denote the point
$\frac{b+\sqrt{N}}{2}$. Then $\bar \p \A \bar \D = \langle
\frac{b+\sqrt{N}}{2} , pa|D| \rangle$, and on these basis the matrix
$B$ of theorem \ref{cm} is
given by $\left(
\begin{array}{cc} 1&0\\ 0 & p \end{array} \right)$. Now $B z_0
= \frac{z_0}{p}$ and $p B^{-1} = \left( \begin{array}{cc} p&0\\
0&1\end{array} \right) = S^{-1}BS$. By theorem \ref{cm},
$g(z_0)^{\sigma(\p)} = [g \circ (pB^{-1})](B z_0)$.

Let $g^\star(z) = g \circ S(z) = g(-1/z)= \theta (\vec{0}, -1/(a|D|z) Q_\B) / \eta(\frac{-1}{|D|z})
\eta(\frac{-1}{z})$. If in (\ref{third}) we replace $z$ by $za|D|$ and $Q_\B$ by
$\text{Adj}\,(Q_\B)$, we get the equation:
\begin{equation}
\theta(\vec{0}, Q_\B (-1/a|D|z)) = (-i)\sqrt{|D|}  az
\theta(\vec{0},\text{Adj}\,(Q_\B) az)
\label{cplx}
\end{equation}

The eta function satisfies the functional equation $\eta(-1/z)=
\sqrt{z/i} \, \eta(z)$. Replacing $z$ by $|D|z$
and multiplying both equations:
$$\eta(-1/z)
\eta(-1/(|D|z)) = \sqrt{|D|} \frac{z}{i} \eta(z) \eta(|D|z)$$
Hence we get:
$$g(-1/z) =  a
\frac{\theta(\vec{0},\text{Adj}\,(Q_\B) az)} {\eta(z) \eta(|D|z)}$$
The $q$-expansion of this quotient has rational coefficients hence it is
fixed by the action of $\sigma_p$, i.e. $g^{\star} \circ \sigma_p  =
g^\star$. Then $[g \circ (pB^{-1})] = g$ and
$(g(z_0))^{\sigma_\p} = g(z_0/p)$.

\begin{prop} With the notation as above, if $\p$ is principal,
$g(z_0)^{\sigma_\p}= g(z_0)$.

\end{prop}

\noindent{\bf Proof.} The proposition reduces to proving that $g(z_0/p)
= g(z_0)$ if $\p$ is principal of norm $p$ which follows from the 
next two lemmas. This completes the proof of theorem \ref{def}
since it implies that $g(z_0)^{\sigma_\p} = g(z_0)$ for all principal
ideals $\p$. $\square$

\begin{lem} \label{theta-action} Let $\bar \p = \langle
\mu \rangle$ be a principal ideal prime to $\A$ and $\bar \D$ of norm $p$ . Then the
theta function $\Theta_\B$ satisfies the formula:
$$\Theta_\B \left(\frac{b+\sqrt{N}}{2ap|D|}\right) = \bar{\mu}
\varepsilon_{\bar \D}(\mu) \left( \frac{p}{|D|} \right) \Theta_\B \left(
\frac{b+\sqrt{N}}{2a|D|}\right)$$
\end{lem}

\begin{note} Since $\varepsilon_{\bar \D}(\mu)
\varepsilon_{\bar \D}(\bar{\mu}) =
\left( \frac{p}{|D|} \right)$, the formula may be written as 
$\Theta_\B (\frac{b+\sqrt{N}}{2ap|D|}) \break = \psi_{\bar \D}(\bar{\mu})
 \Theta_\B(\frac{b+\sqrt{N}}{2a|D|})$
\end{note}

\bigskip

\noindent{\bf Proof.} $\Theta_\B$ is a modular form of weight 1 for
$\Gamma_0(|D|)$ with a quadratic character. We chose $b$ such
that $\bar \p \A \bar \D = \langle
\frac{b+\sqrt{N}}{2}, pa|D| \rangle = \langle \mu
\frac{b+\sqrt{N}}{2}, \mu a |D| \rangle$. Hence there exists a change
of basis matrix $M = 
\left( \begin{array}{cc}
\alpha& \beta\\ 
\gamma & \delta \end{array} \right)$ in $Sl_2(\Z)$ such that $
\left( \begin{array}{cc}
\alpha& \beta\\ 
\gamma & \delta \end{array} \right)
\left( \begin{array}{c}
\frac{b+\sqrt{N}}{2}\\ 
ap|D| \end{array} \right) = 
\left( \begin{array}{c}
\mu \frac{b+\sqrt{N}}{2}\\ 
\mu a |D| \end{array} \right)$.

If $\mu = \frac{m + n \sqrt{N}}{2}$, an easy computation shows that $\delta = \frac{m-nb}{2p}$ and $\gamma
= n|D|a$. In particular $M$ is in $\Gamma_0(|D|)$ and by modularity of
$\Theta_\B$ we have:
$$\Theta_\B \left(\frac{b+\sqrt{N}}{2a|D|}\right) =
\Theta_\B\left(M . \frac{b+\sqrt{N}}{2ap|D|}\right) = \left(\gamma
\frac{b+\sqrt{N}}{2ap|D|} + \delta\right) \chi(\delta)
\Theta_\B\left(\frac{b+\sqrt{N}}{2ap|D|}\right)$$
And the formula:
\begin{equation}
\Theta_\B \left(\frac{b+\sqrt{N}}{2a|D|}\right) = \frac{\mu}{p} \chi(\delta)
\Theta_\B\left(\frac{b+\sqrt{N}}{2ap|D|}\right)
\end{equation}
 where $\chi(d) = \left( \frac{D}{q} \right)$ for any prime $q$ which
is sufficiently large and satisfies $q \equiv d \bmod |D|$. (\cite
{Ogg} Theorem 20, Chapter VI, page 25). Let $q$ be a prime congruent to $1$
modulo $4$ and congruent to $\delta$ modulo $|D|$. Then
$\chi(\delta) = \left( \frac{D}{q} \right) = \left( \frac{|D|}{q}
\right) = \left( \frac{q}{|D|} \right) = \left(
\frac{\frac{m-nb}{2p}}{|D|} \right) = \left(
\frac{\frac{m-nb}{2}}{|D|} \right) \left( \frac{p}{|D|}
\right)$. Then the proof follows from the definition of
$\varepsilon_{\bar \D}$ and the fact that
$\frac{\mu}{p} = (\bar{\mu})^{-1}$. $\square$

\begin{lem} With the same assumptions as above, the eta function
satisfies the equation $\eta(\frac{b+\sqrt{N}}{2p|D|})
\eta(\frac{b+\sqrt{N}}{2p}) = \bar{\mu} \varepsilon_{\bar \D}(\mu) \left(
\frac{p}{|D|} \right) \eta(\frac{b+\sqrt{N}}{2|D|})
\eta(\frac{b+\sqrt{N}}{2})$. 

\noindent In term of ideals:
\begin{equation}
\eta(\bar \p \bar \D) \eta(\bar \p) = \bar{\mu} \varepsilon_{\bar \D}(\mu) \left(
\frac{p}{|D|} \right) \eta(\bar \D) \eta(\scrO _K)
\label{lemma-eta-eq}
\end{equation}
\label{lemma-eta}
\end{lem}
\noindent{\bf Proof.} Since we choose $|N|
\equiv 3 \bmod 4$, and $|N| \not = 3$, the number of units in $H$ is
$2$ (see \cite{Hajir} tables 3 and 4 of page 507). Given a principal
ideal $\langle u \rangle$ with $u \in \scrO _K$, prime to $\langle 6 \rangle$ define:
$$\kappa(u) = \chi_4(N_{K/\Q}(u)) \frac{1}{\bar u}
\frac{\eta^2(u)}{\eta^2(\scrO _K)}$$ 
Where $\chi_4(a) = \left(\frac{-1}{a}\right)$.
Since the number of units in $H$
is $2$, $\kappa$ is a quadratic character (see \cite{Hajir}, Lemma
14). We can write the left hand side of (\ref{lemma-eta-eq}) as:
\begin{equation}\label{eta}
\eta(\bar \p \bar \D) \eta(\bar \p) = \left( \frac{\eta(\bar \p \bar \D)}{\eta(\bar \D)}
\frac{\eta(\scrO _K)}{\eta(\bar \p)} \right)
\frac{\eta^2(\bar \p)}{\eta^2(\scrO _K)} \eta(\scrO _K) \eta(\bar \D)
\end{equation}
If $\mu$ is a generator of $\bar \p$, $\frac{\eta^2(\bar \p)}{\eta^2(\scrO_K)} = \kappa(\mu)
\bar{\mu} \chi_4(p)$. By proposition 10 of \cite{Hajir} 

$\left(\frac{\eta(\bar \p)}{\eta(\scrO _K)} \right)^{\sigma_{ \D}} =
\left( \frac{p}{|D|} \right) \frac{\eta(\bar \p \bar \D )}{\eta(\bar \D)}$. Then we
get:
$$\left( \frac{\eta(\bar \p \bar \D)}{\eta(\bar \D)}
\frac{\eta(\scrO _K)}{\eta(\bar \p)} \right)= \left( \frac{p}{|D|} \right) \left(
\frac{\eta(\bar \p)}{\eta(\scrO _K)} \right) ^{\sigma_{ \D} -1}=\left( \frac{p}{|D|} \right)
\left(\sqrt{\kappa(\mu) \bar \mu \chi_4(p)}\right)^{\sigma_{ \D}-1} $$
By lemma 12 of \cite{Hajir}, $\kappa(-1) = -1$. Since the right term of
(\ref{lemma-eta-eq}) remains unchanged replacing $\mu$ by $-\mu$,
without loss of generality we can choose $\mu$ such that $\kappa(\mu)
= \chi_4(p)$. Replacing each term on the right hand side of
(\ref{eta}) and using Proposition \ref{Artin-Tate} we get:
$$\eta(\bar \p \bar \D) \eta(\bar \p) = \left( \frac{p}{|D|} \right) \varepsilon_{
\D}(\bar \mu) \, \bar \mu \,\eta(\scrO _K) \eta(\bar \D)$$
And the result follows since $\varepsilon_{\D}(\bar \mu) = \varepsilon_{\bar \D}(\mu)$. $\square$

\begin{thm}
The number $n_{\A,\B,\bar \D}$ is in the field $\M_\psi=HT_\psi$. It corresponds to the fields diagram:
\[
\label{Diagram}
\xymatrix{
& \M_\psi \ar@{-}[ld]_h \ar@{-}[rd]^h \ar@{.}[d]\\
H^{} \ar@{-}[d] \ar@{-}[rd]& \M_\psi^{^+} \ar@{.}[ld] \ar@{.}[rd]  &T_\psi \ar@{-}[d] \ar@{-}[ld]\\
 H^{^+} \ar@{-}[rd]_h & K \ar@{-}[d]^2 & T_\psi^{^+} \ar@{-}[ld]^h\\
& \Q
}\]
\begin{center}\text{Diagram 1} \end{center}
\label{field-def}
\end{thm}

\noindent{\bf Proof.} By theorem \ref{def} the number $\theta(\vec{0},z_{\A \bar \D} Q_{\B}) / \eta(z_{\bar \D})
\eta(\scrO_K)$ is in $H$ and $T_\psi$ contains the image of
$\psi_{\bar \D}$ hence $n_{\A,\B, \bar \D}$ is in $\M_\psi$. $\square$ 

\begin{prop} The quotient $\theta_{Q_\B}(z_{\A \bar{\D}})
/ \psi_{\bar\D}(\bar{\A})$ depends only on the class of
$\B$ and the class of $\A$.
\end{prop}

\noindent{\bf Proof.} Independence of $\B$ is clear since
$\Theta_{\B}$ depends only in the class of $\B$. 

\noindent To prove independence of $\A$, let $\alpha \in \scrO_K$ be an element
with prime norm $q$ such that $q \nmid 6a|D|$. By definition $\Theta_\B(z_{\alpha \A \D}) =
\Theta_\B(\frac{b+\sqrt{N}}{2aq|D|})$. Then by lemma \ref{theta-action}:
$$\Theta_\B \left(\frac{b+\sqrt{N}}{2aq|D|}\right) =
\psi_{\bar{\D}}(\bar{\alpha}) \Theta_\B \left(
\frac{b+\sqrt{N}}{2a|D|}\right) \square$$
We will denote by $n_{[\A],[\B],\bar \D}$ the number $n_{\A,\B,\bar \D}$.

\begin{prop} The number $n_{[\A],[\B], \bar \D}$ is an algebraic integer.
\end{prop}

\noindent{\bf Proof.} In theorem \ref{def} we proved that $\theta_{Q_\B}(z_{\A
\bar{\D}})/ (\eta(z_{\bar \D}) \eta(z_{\scrO _K}))$ is an algebraic
integer and the number $\psi_{\bar \D}(\bar \A)$
has norm $N \A$. Since the quotient depends on the class of the ideal
$\A$ but not $\A$ itself, using the Tchebotarev density theorem we can
choose two prime ideals $\p_1$ and $\p_2$ in the same class of $\A$ of
prime norms $p_1$ and $p_2$. Looking at $\p_1$ we see that the minimal
polynomial of $n_{[\p_1],[\B],\bar \D}$ has rational coefficients with only $1$ or $p_1$ in the
denominator. Considering $\p_2$ we see that the minimal polynomial of
$n_{[\p_2],[\B],\bar \D}$ only has $1$ or $p_2$ in the
denominator. Since $n_{[\p_1],[\B],\bar \D} = n_{[\p_2],[\B],\bar \D}$
its minimal polynomial must have integer coefficients.$\square$

\begin{prop} $n_{[\A],[\bar \B],\bar \D}=n_{[\A],[\B],\bar \D}$
\label{conjugationofB}
\end{prop}

\noindent {\bf Proof.} It is easy to check that the theta function
$\Theta_\B$ associated to $\B$ is the same as the theta function
$\Theta_{\text{Adj}\,B}$ associated to the adjoint matrix of
$\B$. Note that $[\B^{-1}]=[\bar \B]$. Clearly the point $z_{\A \bar \D}$ and the number $\psi_{\bar
\D}(\A)$ are independent of $\B$. $\square$

\begin{lem} The character $\psi_{\bar \D}$ satisfy: $\overline{\psi_{\D}(\bar
\A)} = \psi_{\bar \D}(\A)$
\label{character-conjugation}
\end{lem}

\noindent {\bf Proof.} 
$\overline{\psi_{\D}(\bar\A)} \psi_{\D}(\bar\A) = N\A$, and
 $N\A = \psi_{\D}(\bar\A) \psi_{\D}(\A)
\varepsilon_{\D}(N\A)$ hence $\overline{\psi_{\D}(\bar\A)} =
\left( \frac{N\A}{|D|}\right) \psi_{\D}(A)$. We chose the
characters so that $\psi_{\bar \D}(\A) = \psi_{\D}(\A) \varepsilon_{\bar \D}(\A^h)
\varepsilon_{\D}(\A^h) \break = \psi_\D(\A) \left( \frac{(N\A)^h}{|D|}
\right)$ (see (\ref{consistency})). Since $|N|$ is prime, $h$ is odd. $\square$

\begin{prop} $\overline{n_{[\A],[\B],\bar \D}}=n_{[\bar \A],[\B],\D}$
\label{conjugacy-nonppal}
\end{prop}

\noindent {\bf Proof.} It is clear from their definition that
$\overline{\Theta_{\B}(z_{\A \bar \D})} = \Theta_{\B}(-\overline{z_{\A
\bar \D}})$ and $\overline{\eta(z_{\A \bar \D})} \break =
\eta(-\overline{z_{\A \bar \D}})$ . Since $-\overline{z_{\A \bar \D}}
= z_{\bar \A \D}$ and $\overline{\psi_{\D}(\bar \A)} = \psi_{\bar
\D}(\A)$, the result follows. $\square$.

\begin{prop} If the ideal $\D$ is principal in $\scrO_K$,
$\overline{n_{[\A],[\B],\bar \D}} = n_{[\bar \A],[\B],\bar \D}$
\label{conjugacy}
\end{prop}

\noindent{\bf Proof.} The proof of this proposition involves the same
kind of techniques used on the previous ones (a little more tedious) so we
omit the proof. $\square$

In particular this implies that if $\A$ and $\D$ are both principal
then the number $n_{[\A],[\B],\bar \D}$ lives in a subfield of $\M_\psi$
which we denote $\M_\psi^+$ (following \cite{Buhler-Gross} notation, see
page 13) and corresponds to the previous field diagram (see theorem
\ref{Diagram} Diagram 3.3). We will be needing the next lemmas for the
theorem relating the numbers $n_{[\A],[\B],\bar \D}$ for different
ideals $\D$.

\begin{lem} Let $\D$ and $\D'$ be two prime ideals of $\Q(\sqrt{N})$
with norm $|D|$ and $|D'|$ respectively, and let $\mu \in
\Q(\sqrt{N})$ be such that $\mu \D = \D'$. Then $\frac{\eta^2(\A \D')}{\eta^2(\A \D)} =
\bar \mu \kappa(\mu) \chi_4(N\mu)$.
\label{etas-squared}
\end{lem}

\noindent{\bf Proof.} Note that although $\kappa$ is defined on
integer elements, since it is a character on $\left(\scrO_K / 12
\scrO_K \right) ^ \times$, we can extend it multiplicatively to all
elements in $\Q(\sqrt{N})$ with both numerator and denominator prime
to $12$. By definition $\kappa(\mu)
= \frac{1}{\bar \mu} \chi_4(N\mu)
\frac{\eta^2(\mu)}{\eta^2(\scrO_K)}$ then we are led to prove that
$\frac{\eta^2(\A \D')}{\eta^2(\A \D)}
\frac{\eta^2(\scrO_K)}{\eta^2(\mu)} = 1 $.\\
By Proposition 10 of \cite{Hajir} we can write the left hand side as
$\left( \frac{\eta^2(\A \D)}{\eta^2(\scrO_K)} \right)
^{\sigma_{(\bar{\D}' \bar{\D}^{-1})} -1}$.

Since $\frac{\eta^2(\A \D)}{\eta^2(\scrO_K)}$ is in $H$ (by
theorem 20 of \cite{Hajir}) then $\sigma_\A$ represents the classical
Artin-Frobenius map from $Cl(\scrO_K)$ to $\Gal(H/K)$, and since $\bar{\D}'
\bar{\D}^{-1}$ is principal, $\sigma_{\bar{\D}' \bar{\D}^{-1}}$ is the
identity. $\square$

\begin{lem} Let $\D$ and $\D'$ be two prime ideals of $\Q(\sqrt{N})$
 such that $\D \sim \D'$. Then
$\frac{\eta(\A \D') \eta(\D)}{\eta(\A \D) \eta(\D')} =
\varepsilon_{\D}(\bar \A^h) \varepsilon_{\D'}(\bar \A^h)$
\label{more-etas}
\end{lem}

\noindent{\bf Proof.} By proposition 10 of \cite{Hajir} we have:
\begin{equation}
\frac{\eta(\A \D') \eta(\D)}{\eta(\D') \eta(\A \D) } = \left(
\frac{\eta(\A)}{\eta(\scrO_K)} \right)^{\sigma_{\bar \D'}} \left(
\frac{\eta(\A)}{\eta(\scrO_K)} \right)^{-\sigma_{\bar \D}}
\left(\frac{a}{|D|}\right) \left( \frac{a}{|D'|} \right)
\label{pf-etas}
\end{equation}

Since the Artin-Frobenius map is a homomorphism:
$$\left( \frac{\eta(\A)}{\eta(\scrO_K)} \right)^{\sigma_{\bar
\D'}-\sigma_{\bar \D}} = \left( \left( \frac{\eta(\A)}{\eta(\scrO_K)}
\right)^{\sigma_{(\bar \D'(\bar\D)^{-1})}-1}\right)^{\sigma_{\bar
\D}}$$ But $\left( \frac{\eta(\A)}{\eta(\scrO_K)}
\right)^{\sigma_{(\bar \D'(\bar\D)^{-1})}-1} = \pm1$ (see the proof of
lemma \ref{etas-squared}), then $\sigma_{\bar \D}$ acts trivially on
it.  Let $\mu \in \Q(\sqrt{N})$ be such that $\bar \D' \bar \D^{-1}$
is the principal ideal generated by $\frac{\bar \mu}{|D|}$ then by theorem 19 of \cite{Hajir}:
$$\left(\frac{\eta(\A)}{\eta(\scrO_K)} \right)^{\sigma_{(\bar \D'(\bar
{\D})^{-1})}-1} = \kappa\left(\frac{ \mu}{|D|}\right)^{\frac{a-1}{2}}
\left( \frac{\bar \mu |D|}{\bar \A} \right) $$
Since $|D|$ is prime to $12$, and $\kappa$ is a multiplicative
quadratic character, $\kappa(\frac{\mu}{|D|})=\kappa(\mu)
\kappa(|D|)$. The character $\kappa$ defined on $\left(\scrO_K/12
\scrO_K \right)^\times$ factors as a product of two characters,
$\kappa_3$ from $\left(\scrO_K/3 \scrO_K \right)^\times$ to the group
of third roots of unity and $\kappa_4$ from $\left(\scrO_K/4 \scrO_K
\right)^\times$ to the group of fourth roots of unity (see lemma 14 of
\cite{Hajir}). In our
case $\kappa_3 =1$ and the character is completely determined from the
congruence modulo $4$. Then $\kappa(|D|) = \kappa(-1) =-1$. Using
the quadratic reciprocity law,
\begin{equation}
\left(\frac{\eta(\A)}{\eta(\scrO_K)} \right)^{\sigma_{(\bar \D'(\bar
{\D})^{-1})}-1} = \kappa(
\mu)^{\frac{a-1}{2}} \left( \frac{\bar \mu}{\bar \A} \right) \left(
\frac{a}{|D|} \right)
\end{equation}
Also since $\kappa(\mu) \kappa(\bar \mu) =
\kappa(|D||D'|) = 1$, $\kappa(\mu) = \kappa(\bar \mu)$ and we can
write (\ref{pf-etas}) as:
$$\frac{\eta(\A \D') \eta(\D)}{\eta(\A \D)\eta(\D')} = \kappa(\bar
\mu)^{\frac{a-1}{2}} \left( \frac{\mu}{\A} \right) \left( \frac{a}{|D'|}
\right) $$
Since $\bar \D \D'$ is the principal ideal generated by $\mu$ and $\varepsilon$ is a multiplicative
quadratic character,
\begin{equation}
\varepsilon_{\D}(\bar \A ^h)
\varepsilon_{\D'} (\bar \A^h) = \varepsilon_{\D}(\bar \A^h)
\varepsilon_{\bar \D}(\bar \A^h) \varepsilon_{\bar \D \D'}(\bar \A^h) 
= \left( \frac{a}{|D|} \right) \left( \frac{\bar \A^h}{\mu} \right)
\end{equation}
Using the reciprocity law in $\Q(\sqrt{N})$ (see
for example theorem 21 of \cite{Hajir}):
\begin{equation}
 \left( \frac{\bar \A^h}{\mu}\right) = \left( \frac{\mu}{\bar \A^h}
 \right) \kappa(\bar \mu)^{\frac{a-1}{2}} =
\left( \frac{\mu}{\bar \A}
 \right) \kappa(\bar \mu)^{\frac{a-1}{2}}=
 \kappa(\bar \mu)^{\frac{a-1}{2}} \left( \frac{\mu}{\A}\right) 
  \left( \frac{|D||D'|}{a} \right) 
\end{equation}
And the lemma follows from $\left( \frac{|D||D'|}{a} \right) = \left(
 \frac{a}{|D|}\right) \left(\frac{a}{|D'|} \right)$. $\square$

\begin{lem} Let $A: \R^{2n}\times \R^{2n} \rightarrow \R$ be the
 skew-symmetric form given by the matrix $A := \left(
\begin{array}{cc}
0 & I_n\\ -I_n & 0 \end{array} \right)$. Then the following data on $\R^{2n}$
 are equivalent:
\label{Mumfordeq}
\begin{enumerate}
\item a complex structure $U:\R^{2n} \rightarrow \R^{2n}$ (i.e. a
  linear map with $U^2=-I_n$) such that there exists a positive
  definite Hermitian form $H$ for this complex structure with
  imaginary part $A$.
\item an $n$-dimensional complex subspace of $\C^{2n}$ such that if we
  note $A_\C$ the complex linear extension of $A$, we have:
\begin{itemize}
\item $A_\C(x,y) = 0$ for all $x,y$ in the subspace.
\item $i A_\C(x,\bar x)<0$ for all nonzero $x$ in the subspace.
\end{itemize}
\item a complex matrix $\Omega$ in $\Siegel_n$
\end{enumerate}
\end{lem}
This are three of the four equivalent conditions proved on Lemma 4.1 of
\cite{Mumford}. The equivalence associates to $\Omega \in \Siegel_n$
the image of the map $X \mapsto (X,- \Omega X)$ as an $n$-dimensional
subspace of $\C^{2n}$.

\begin{thm} Let $z_{\A \D} Q_\B$ and $z_{\A \D'} Q_{\B'}$ be two points
in $\Siegel_2$ such that they are equivalent modulo $\Gamma_{12}$ and $\D
\sim \D'$ in $\Q(\sqrt{N})$. Then $n_{[\A],[\B],\bar D} = \pm n_{[\A],[\B'],\bar D'}$
\label{comparison}
\end{thm}

\noindent {\bf Proof.} For simplicity we will
denote $\Omega_{\D} := z_{\A \D} Q_{\B}$ and $\Omega_{\D'} := z_{\A \D'}Q_{\B'}$.
Since $\Omega_{\D}$ is
equivalent to $\Omega_{\D'}$ there exists a matrix $ \gamma =\left(
\begin{array}{cc}
A&B\\ C&D \end{array} \right)$ in $Sp_4(\Z)$ such that $\gamma \star
(\Omega_\D) = \Omega_{\D'}$. By the previous lemma, giving a point
$\Omega_{\D}$ in the Siegel space is equivalent to giving the subspace
of $\C^4$ $(I_2,-\Omega_\D)^t$ where the action of $Sp_4(\Z)$ is given
by multiplication on the left by $(\gamma^t)^{-1}$. Then $\gamma \star
\left( \frac{I_2}{-\Omega_\D}\right) = \left(
\begin{array}{cc} D&-C\\ -B&A \end{array} \right)
\left(\frac{I_2}{-\Omega_\D}\right)= \left( \frac{ C
\Omega_\D+D}{-(A\Omega_\D+B)} \right) = \left(
\frac{I_2}{-\Omega_{\D'}} \right) (C \Omega _\D +D)$\\
By the coherent way we chose characters,
$\frac{\psi_{\D}(\A)}{\psi_{\D'}(\A)} = \varepsilon_\D(\A^h)
\varepsilon_{\D'}(\A^h)$. Hence:
$$\frac{n_{[\A],[\B],\bar D}}{ n_{[\A],[\B'],\bar D'}} =
\frac{\theta(\Omega_\D)}{\theta(\Omega_{\D'})}
\frac{\eta(\D')}{\eta(\D)} \varepsilon_{\D}(\bar \A^h)
\varepsilon_{\D'}(\bar \A^h)=
\frac{\theta(\Omega_\D)}{\theta(\Omega_{\D'})}
\frac{\eta(\A \D')}{\eta(\A \D)}$$
The last equality follows from lemma \ref{more-etas}. We
claim that:
\begin{equation}
\frac{\theta^2(\Omega_\D)}{\theta^2(\Omega_{\D'})} = \Det(C \Omega_\D
+D)^{-1} = \frac{\eta^2(\A \D)}{\eta^2(\A \D')}
\end{equation}
The first equality follows at once from the functional equation of the
theta function. Since $|D|$ is prime and $\Det(Q) = |D|$ there exists
matrices $U,V \in Sl_2(\Z)$ such that $UQV = \left( \begin{array}{cc}
1&0\\ 0&|D| \end{array} \right)$ (respectively $U'$ and $V'$ for
$Q'$). Then:
$$
\left( \begin{array}{cc}
V^{-1}&0\\
0&U \end{array} \right) 
\left( \begin{array}{c}
I_2\\
-\Omega_\D \end{array} \right) V =
\left( \begin{array}{c}
I_2\\
-UQV z_{\A\D} \end{array} \right) = 
\left( \begin{array}{cc}
1&0\\
0&1\\
-z_{\A \D}&0\\
0 & -z_\A \end{array} \right) 
$$
Similarly:
$$
\left( \begin{array}{cc}
V'^{-1}&0\\
0&U' \end{array} \right) 
\left( \begin{array}{c}
I_2\\
-\Omega_{\D'} \end{array} \right) V' =
\left( \begin{array}{c}
I_2\\
-U'Q'V' z_{\A\D'} \end{array} \right) = 
\left( \begin{array}{cc}
1&0\\
0&1\\
-z_{\A \D'}&0\\
0 & -z_\A \end{array} \right) 
$$
We split in two cases:
\begin{itemize}
\item if $\D' = \bar \D$ we take basis $\D =
\langle |D|, \frac{b+\sqrt{N}}{2} \rangle$ and $\A = \langle a,
\frac{b+\sqrt{N}}{2} \rangle$. Let $r$ be such that $r |D| \equiv b
\bmod a$ then $\D' = \langle |D'|, \frac{(2r|D|-b)+\sqrt{N}}{2}
\rangle$ and $\A \D' = \langle a|D'|, \frac{(2r|D|-b)+\sqrt{N}}{2}
\rangle$. Let $\mu \in K$ be such that $\mu \D = \D'$,
then $\A \D' = \langle a|D'|,\frac{(2r|D|-b)+\sqrt{N}}{2} \rangle =
\langle \mu a|D|, 
\mu (\frac{b+\sqrt{N}}{2}) \rangle = \mu \A \D$ hence there exists a
matrix $M = 
\left( \begin{array}{cc}
\alpha&\beta\\
\gamma& \delta \end{array} \right)$ in $Sl_2(\Z)$ such that:
$
M 
\left( \begin{array}{c}
\mu(\frac{b+\sqrt{N}}{2})\\
\mu a|D| \end{array} \right)=
\left( \begin{array}{c}
\frac{(2r|D|-b)+\sqrt{N}}{2}\\
a|D'| \end{array} \right)
$

\item if $\D' \neq \bar \D$, we may choose basis 
$\D = \langle |D|, \frac{b+\sqrt{N}}{2} \rangle$, $\D' = \langle |D'|,
\frac{b+\sqrt{N}}{2} \rangle$ and $\A = \langle a,
\frac{b+\sqrt{N}}{2} \rangle$. If $\mu$ is such that $\mu \D = \D'$,
then $\A \D' = \langle a|D'|,\frac{b+\sqrt{N}}{2} \rangle = \langle \mu a|D|,
\mu (\frac{b+\sqrt{N}}{2}) \rangle = \mu \A \D$ hence there exists a matrix $M =
\left( \begin{array}{cc}
\alpha&\beta\\
\gamma& \delta \end{array} \right)$ in $Sl_2(\Z)$ such that:
$
M 
\left( \begin{array}{c}
\mu(\frac{b+\sqrt{N}}{2})\\
\mu a|D| \end{array} \right)=
\left( \begin{array}{c}
\frac{b+\sqrt{N}}{2}\\
a|D'| \end{array} \right)
$
\end{itemize}
In both cases, let $N := \left( \begin{array}{cccc}
\delta&0&-\gamma&0\\
0&1&0&0\\
-\beta&0&\alpha&0\\
0&0&0&1 \end{array} \right)$, then it is clear that:
$$
N
\left( \begin{array}{cc} 
1&0\\
0&1\\
-z_{\A \D}&0\\
0 & -z_\A \end{array} \right) 
\left( \begin{array}{cc}
\frac{\mu|D|}{|D'|}&0\\
0&1 \end{array} \right) =
\left( \begin{array}{cc}
1&0\\
0&1\\
-z_{\A \D'}&0\\
0 & -z_\A \end{array} \right) 
$$
Combining these results we get that:
$$
\left( \begin{array}{cc}
V'&0\\
0&U'^{-1} \end{array} \right)
N
\left( \begin{array}{cc}
V^{-1}&0\\
0&U \end{array} \right) 
\left( \begin{array}{c}
I_2\\
-\Omega_\D \end{array} \right) V 
\left( \begin{array}{cc}
\frac{\mu|D|}{|D'|}&0\\
0&1 \end{array} \right)
V'^{-1} = 
\left( \begin{array}{c}
I_2\\
-\Omega_{\D'} \end{array} \right)
$$
and
$$
\left( \begin{array}{c}
I_2\\
-\Omega_{\D'} \end{array} \right)
= \left( \begin{array}{cc}
D&-C\\
-B&A \end{array} \right)
\left( \begin{array}{c}
I_2\\
-\Omega_\D \end{array} \right)
(C\Omega_\D+D)^{-1}
$$
Since both lattices have the same volume then $|\Det(C
\Omega_\D+D)|^{-1} = \frac{|\mu||D|}{|D'|}$.\\
By lemma \ref{etas-squared}, $\frac{\eta^2(\A
\D)}{\eta^2(\A \D')} = \frac{1}{\bar \mu} \kappa(\mu)=
\frac{\mu|D|}{|D'|} \kappa(\mu)$. Now $\Det(C \Omega_D+D)^{-1}$
and $\kappa(\mu) \frac{\mu|D|}{|D'|}$ have the same absolute value and
both lie in $\Q(\sqrt{N})$ hence they differ by $\pm1$. Then 
$$\left(
\frac{\theta(\Omega_\D)}{\theta(\Omega_{\D'})} \frac{\eta(\A
\D')}{\eta(\A \D)} \right)^2 = \Det(C \Omega_\D +D)^{-1} \bar \mu
\kappa(\mu) = \pm1$$
Taking square roots:
$$
\sqrt{\pm 1} = \frac{\theta(\Omega_\D)}{\theta(\Omega_{\D'})}
\frac{\eta(\A \D')}{\eta(\A \D)}
$$
By theorem \ref{field-def} we know
that $\frac{\theta(\Omega_\D)}{\eta(\D) \eta(\scrO_K)}$ and
$\frac{\theta(\Omega_\D')}{\eta(\D') \eta(\scrO_K)}$ are in $H$.
Since $\sqrt{-1} \not \in  H$ the theorem follows. $\square$
\medskip

\noindent It is not clear how to determine the sign a priori, and we are not
able to give any answer in this direction.
\section{Equivalence of special points}

The problem of determining whether two points in $\Siegel_2$ are
equivalent or not is complicated in general. For our case we will get
this equivalence via ideals in quaternion algebras. A good reference
for the basic definitions and some elementary facts about quaternion
algebras is Pizer's paper (\cite{Pizer}).

Let $B$ be a quaternion algebra over $\Q$. A lattice $\lat$ is a rank $4$
$\Z$-module. An order $O$ is a lattice that is a ring with
unity. Given an order $O$ a left $O$-ideal is a lattice $\lat$ such that
$\lat_p:= \lat \otimes _\Z \Z_p = O_p \alpha_p$ where $\alpha_p$ is an
element in $B_p^\times$. Given a lattice $\lat$ we define its left
 order $O_l(\lat):=\{x \in B \, | \, x\lat \subset
\lat\}$ (respectively the right order). We define $\Norm(\lat)$ as the positive generator of the $\Z$-module
$\langle \Norm(x) \, | \, x \in \lat \rangle$

\begin{prop} 
Let $B$ be a quaternion algebra over $\Q$ ramified at $p_1 , \ldots,
p_n$ and $\lat$ be an ideal in
$B$. Then $O_l(\lat)$ is a maximal order if and only if disc($\lat$)$=(p_1
\ldots p_n)^2 \Norm(\lat)^4$
\label{order-level}
\end{prop}

\noindent{\bf Proof.} By definition disc($\lat$) is the determinant of
the bilinear form associated to $\lat$ on any basis. Since $\lat$ is locally
principal at all primes, given a finite prime $q$, $\lat_q = O_l(\lat)_q
\alpha_q$. Clearly disc$(\lat_q) = \Norm(\alpha_q)^4$ disc$(O_q)$; then the
statement follows from the fact that this proposition is true
replacing $\lat$ by an order $O$ and $\Norm(\lat)$ by $1$ (see \cite{Pizer}
Proposition 1.1, page 344), and the fact that the norm of $\lat$ is the
product over all primes $q$ of $q^{v_q(N\alpha_q)}$ where $v_q(n)$ is
the $q$-valuation. $\square$

We restrict ourselves to the case $B$ a quaternion algebra over $\Q$
ramified at the prime $|N|$ and infinity.

\begin{lem} Let $O$ be a maximal order, $\{I_1,
\ldots , I_h\}$ a set of left $O$-ideal representatives, and $\{R_1,
\ldots , R_h\}$ be the right orders of $\{I_1 , \ldots
, I_h\}$ respectively. Then for a given $i=1, \ldots , h$ the maximal
order $R_i$ appears twice on the list if and only if there is no
embedding of $\Z[\sqrt{N}]$ into $R_i$.
\label{Eichler-lemma}
\end{lem}

\noindent{\bf Proof.} Although this is a well known statement we give
a proof since we will use it latter. An embedding of $\Z[\sqrt{N}]$
into $R_i$ is determined by the image of $\sqrt{N}$. Hence giving such
an embedding is equivalent to giving an element $\beta \in R_i$ of trace
zero and norm $|N|$.
Let $\bfP$ be the bilateral $O$-ideal of norm
$|N|$ . For a given left $O$-ideal $I_j$, the ideal $\bfP I_j$ is
another left $O$-ideal. Note that if $\bfP_j$ is the bilateral $R_j$ ideal of norm
$|N|$, then $I_j^{-1} \bfP I_j = \bfP_j$ by the uniqueness of
such a bilateral ideal. Then the ideals $I_j$ and
$\bfP I_j$ are equivalent if and only if there exists $\beta \in
R_j^\times$ such that $I_j \beta = \bfP I_j$. Multiplying on the left
by $I_j^{-1}$ we see that $R_j \beta = I_j^{-1} \bfP I_j = \bfP_j$
hence $\bfP_j$ is principal, and the element $\beta$ has norm
$|N|$. Since $|N|$ is a ramified prime i.e. $B_{|N|}$ is a division
ring, it is easy to see that if $\Norm(\alpha) = |N|$ then $\Tr(\alpha) =
0$.  

To see that this is the only way in which a maximal order $R$ appears
twice on the list of right orders, suppose that $I$ and $J$ are two
nonequivalent left $O$-ideals with same right order $R$. Then
$I^{-1}J$ is a non-principal bilateral ideal for $R$. Let $\bfP_R$ be
the ideal of norm $|N|$ in $R$, then $\bfP_R$ is non-principal and $J$ is equivalent to
$\bfP I$. $\square$

\subsection{Siegel Space and applications}

\begin{defn} Let $\lat$ be a $\Z$ lattice of rank $2n$ and $V$ the vector
  space $\lat \otimes \R$. We call a triple $(P,J,U)$ a Siegel point if:
\begin{itemize}
\item $P$ is a real $2n \times 2n$ symmetric matrix such that the
  associated quadratic form $P(x,y)$ is positive definite (that will correspond to the real part of $H$ ).
\item $J$ is a real $2n \times 2n$ non-degenerate skew symmetric
  matrix with associated form $J(x,y)$
(that will correspond to the imaginary part of $H$).
\item $U \in \R^{2n \times 2n}$ is such that $U^2 = -I_{2g}$. (complex structure)
\end{itemize}
with the relation:
\begin{equation}
-JU = U^tJ=P
\label{hermitian-condition}
\end{equation}
\end{defn}
Via the matrix $U$ we can put a complex structure on the vector space
$V$. Let $H$ be the bilinear form $H(x,y):= P(x,y)+iJ(x,y)$. The
condition (\ref{hermitian-condition}) implies that $H(ix,y) =
iH(x,y)$.  Since $J$ is skew symmetric and $P$ symmetric, it follows
that $H(x,y) = \overline{H(y,x)}$.  Then $H$ defined in this way is a
positive definite Hermitian form. Each choice of a reduced basis for
$J$ will give a point in the Siegel space (by Lemma \ref{Mumfordeq})
and different bases give equivalent points.

Given two lattices $\lat$ and $\lat'$, a {\it morphism} $\gamma:\lat \rightarrow
\lat'$ is an $\Z$-linear map from $\lat$ to $\lat'$. Given
$\gamma:\lat' \rightarrow \lat$ an isomorphism of lattices, we define an
action of $\gamma$ on a Siegel point $(P,J,U)$ as
$(\gamma^*P,\gamma^*J,\gamma^*U)$ where given $x,y \in \lat'$,
$\gamma^*P(x,y) = P(\gamma(x),\gamma(y))$, $\gamma^*J(x,y) =
J(\gamma(x),\gamma(y))$ and $\gamma^*(x)=\gamma^{-1}(U(\gamma(x)))$.\\
If we choose $V_0$ to be a skew
symmetric reduced base for $J$, i.e. a base where $J$ is of the form
$\left(
\begin{array}{cc} 0&I_n\\ -I_n&0
\end{array} \right)$, and $\gamma$ is an automorphism sending a skew
symmetric reduced basis to another one, then $\gamma \in Sp_{2n}(\Z)$ and the action
of $\gamma$ on the Siegel point $\Omega$ associated to $V_0$ is
the usual action of $Sp_{2n}(\Z)$ on $\Siegel_n$.

\subsection{Siegel Points from Quaternion algebras} 

Let $N$ be the negative of a prime congruent to $3$ modulo $4$, and $B = (-1,N)$  the
quaternion algebra ramified at $N$ and infinity. Let $O$ be a maximal
order in $B$ such that there exists an embedding (not necessarily
optimal) of $\Z + \Z \sqrt{N}$ into $O$. Let $u \in O$ be the
image of $\sqrt{N}$, i.e. $u^2 = N$ and $\Tr(u) = 0$. By $I$ we will
denote a left $O$-ideal for a maximal order $O$. To $I$ we associate a
Siegel point $(P,J,U)_I$ as follows: 

\noindent $\bullet$ We take $V$ the real vector space $V:= B \otimes_\Q
\R$.\\
$\bullet$ Define $U$ acting on $V$ as left multiplication by 
$\frac{u}{\sqrt{|N|}}$.\\
$\bullet$ We think of $I$ as a full rank lattice in $V$.\\
$\bullet$ For $x,y \in I$ define $P(x,y) := \frac{1}{\sqrt{|N|}} \Tr(x \bar
y)/\Norm(I)$.\\
$\bullet$ For $x,y \in I$ define $J(x,y) := \Tr(u^{-1} x \bar y)/\Norm(I)$.

\begin{prop} The triple $(P,J,U)_I$ defined as above is a Siegel point.
\end{prop}

\noindent {\bf Proof.} We start checking the properties of the
matrices $P$, $J$ and $U$:

\noindent $\bullet$ $P$ is a real form. Since $\Tr(x \bar y)$ is real,
$\Tr(x \bar y) = \Tr(y \bar x)$ which implies that $P(x,y)$ is
symmetric. Clearly $P(x,x) = \frac{1}{\sqrt{|N|}} \Norm(x)/\Norm(I)$ is
positive definite.\\
$\bullet$ $J$ is a real form. Since $u$ is pure imaginary, $u^{-1}$
is also. Then $J(x,x) = \Tr(u^{-1}\Norm(x))/\Norm(I)=0$. It is also clear that
$J(x,y)$ is non-degenerate, since for any non-zero $x \in V$,
$J(x,u^{-1}x) \neq 0$. Since $J(x,x)=0$ for all $x$ it follows that
$J(x,y)=-J(y,x)$.\\
$\bullet$ Let $x \in V$, then $U^2(x) = U(\frac{u}{\sqrt{|N|}} x) =
\frac{u^2}{|N|} x = -x$.\\ 
As for the relation, it is easy to check that $J(\frac{u}{\sqrt{|N|}} x,y) =
P(x,y)$ and that $J(x,\frac{u}{\sqrt{|N|}}y) = -P(x,y)$. $\square$
\begin{defn}
Given a lattice $\lat$ in $B$ we define its dual by $\lat^\# := \{ b \in B :
\Tr(b\lat) \subset \Z \}$. Given an order $R$ we define its different by
$R^{\iota} := N R^\#$. 
\end{defn}
\begin{prop}
If $O$ is a maximal order, $O^\iota$ is a bilateral
ideal for $O$ of index $N^2$, and $\frac{1}{N} O \subset O^{\iota}
\subset O$.
\label{different}
\end{prop}
\noindent{\bf Proof.} See \cite{Vigneras} Lemma 4.7, page 24.

\begin{prop} If $x,y \in I$ then $J(x,y) \in \Z$. Also the matrix
of $J$ on the basis given by $I$ has determinant $1$.
\end{prop}
\noindent {\bf Proof.} Since we are considering the reduced norm, if
$V$ is the matrix associated to multiplication (on the left or on the
right) by $v$, then $\Norm(v) = \sqrt{\det(V)}$. Let $W(x,y):= \Tr(x \bar
y)$ be the bilinear form of $B$. If we denote $W$ the matrix of
$W(x,y)$ on the basis given by $I$, $J = \frac{1}{\Norm(I)} (U^{-1})^t W$.
Then $\det(J) = \Norm(I)^{-4}\Norm(u)^{-2} \det(W)$. By definition $\det(W) =
disc(I)$, which is an ideal for a maximal order, then by Proposition
\ref{order-level} $disc(I) = N^2 \Norm(I)^4$ and $\det(J)=1$.

Since the trace is linear, $J(x,y) = \Tr(u^{-1} x \frac{\bar
y}{\Norm(I)})$. For ideals $I$ with maximal left order it is true that
$I^{-1} = \bar I /\Norm(I)$ and $I I^{-1} = O$, hence $J(x,y) \in \Z$ for
all $x,y \in I$ if and only if $\Tr(u^{-1} v) \in \Z$ for all $v \in
O$. By proposition \ref{different} this is the same as $u^{-1} \in
O^\#$. But $u^{-1} = - \frac{u}{N}$, and since $u \in O$ it follows
that $\frac{u}{N} \in \frac{1}{N}O \subset O^\#$. $\square$

This gives a method for assigning to every left $O$-ideal a Siegel
point. Note that choosing different skew symmetric reduced basis of
$I$ will give equivalent Siegel points. From now on we fixed a maximal
order $O$ with an embedding of $\Z[\sqrt{N}]$.

\begin{prop} Let $u \in O$ with $\Norm(u)=|N|$ and $\Tr(u)=0$, and denote
  by $U$ the complex multiplication associated to $u$. 
  If $I$, $I'$ are two equivalent left $O$-ideals, then
  the Siegel points $(P,J,U)_I$ and $(P,J,U)_{I'}$ are equivalent.
\label{equivalence}
\end{prop}
\noindent {\bf Proof.} Since $I \sim I'$ there exists $\alpha \in
B^\times$ such that $I=I' \alpha$. Let $W$ denote the isomorphism of
$B$ given by $W(v) = v\alpha$. We claim that $W$ is the isomorphism
that makes the two Siegel points equivalent.\\ Since $W(I') = I$, we
need to check that $W^* P = P'$, $W^* J =J'$ and $W^*U=U$.

$\bullet$ If $x,y \in I'$ by definition $(W^*P)(x,y) := P(W(x),W(y))
= P(x \alpha, y \alpha) = \frac{\Tr(x \alpha \bar{\alpha} \bar y)}{\Norm(I)}
= \frac{N \alpha}{\Norm(I)} \Tr(x \bar y) = P'(x,y)$.

$\bullet$ The equality $W^*J=J'$ follows from a similar argument.

$\bullet$ By definition $U$ is given by multiplying on the left by
  $u/\sqrt{|N|}$ while $W$ is given by multiplying on the right by
  $\alpha$ then clearly this maps commute with each other and
  $W^*U:=W^{-1}\circ U \circ W = U$. $\square$
\begin{lem} Let $U$ be the complex multiplication given by $u$
and $\alpha \in B$ an element such that $\alpha O \alpha^{-1} = O$. Define $I'=
\alpha I \alpha^{-1}$ and $u' = \alpha u \alpha^{-1}$, then $(P,J,U)_I
\sim (P',J',U')_{I'}$.
\label{equiv}
\end{lem}
\noindent{\bf Proof.} Let $W:B \rightarrow B$ be the isomorphism
defined by $W(x) = \alpha x \alpha^{-1}$. By hypothesis $W(R) = R$,
$W(I)=I'$. It is easy to see that $W^*P = P'$ and $W^* J = J'$.
If $x \in B$ then $W^{-1}\circ U
\circ W(x)= W^{-1} \circ U(\alpha x \alpha^{-1}) = W^{-1}(u \alpha x
\alpha^{-1})/\sqrt{|N|} = \alpha^{-1} u \alpha x/\sqrt{|N|} =
U'(x)$. $\square$

This lemma suggests that we should consider not just elements $u$ in
$O$ corresponding to $\sqrt{N}$ (i.e. $u^2=N$ and $\Tr(u)=0$) but
modulo conjugation by the normalizer of $O$. It is clear that $Norm(O) =
\{ h \in B \, | \, Oh \text{ is bilateral}\}$. All bilateral ideals are
principal, generated by $u^sm$ where $s=0,1$ and $m$ is a rational
number (see \cite{Eichler2} proposition 1, page 92). The generator of
an ideal is well defined up to units in $O$, then $Norm(O) = \{\zeta u^s
m \, | \, s=0\text{ or }1\, , m \in \Q \text{ and } \zeta \in O \text{
is a unit}\}$
\begin{cor}
If $I$ and $I'$ are left $O$-ideals with the same right order then the
Siegel points $(P,J,U)_I$ and $(P,J,U)_{I'}$ are equivalent.
\label{sameRorder}
\end{cor}

\noindent{\bf Proof.} If $I$ and $I'$ are equivalent this follows from
proposition \ref{equivalence}. If $I$ and $I'$ are not equivalent, we
know by lemma \ref{Eichler-lemma} that $O_r(I)$ has no embedding
of $\Z[\sqrt{N}]$. Let $u$ be the element in $O$ giving the complex
multiplication. Then $uI$ has the same left and right order as $I$ but
they are not equivalent, hence $uI \sim I' \sim uIu^{-1}$. By
proposition \ref{equiv} the Siegel points $(P,J,U)_I$ and
$(P,J,U')_{uI}$ are equivalent. Just note that $U'$ is given by $u^{-1}u u =
u$. $\square$

In particular we should index the Siegel points not by the class
number of ideals, but by the type number of maximal orders.
We still have equivalent Siegel points coming from conjugation by
units of $O$ and this are all the possibilities for $Norm(O)$. For
counting equivalent classes of Siegel points, fixed a maximal order
$O$ we have to count the number of embeddings of $\Z[\sqrt{N}]$ into
$O$ modulo conjugation by units of $O$.

Given a maximal ideal $O$, let $\Ba:=\{I_1 , \ldots , I_h\}$ be a set of
left $O$-ideal representatives and $\Ta:=\{ R_1 , \ldots, R_t\}$ the
distinct right orders of the ideals in $\Ba$. We index the Siegel points
by pairs $(\phi,R_i)$ where $\phi$ is an embedding from $\Z[\sqrt{N}]$
to some $R_j$ and $R_i$ is an order in $\Ta$. By this we mean the Siegel
point obtained with the complex multiplication given by
$\phi(\sqrt{N})$, and an ideal $I$ with left order $R_j$ and right
order $R_i$.

If $d$ is a negative discriminant we denote by $h(d)$ the class number
of binary quadratic forms of discriminant $d$. Let $u(d) =1$ unless
$d=-3,-4$ when $u(d)=3,2$ respectively (half the number of units in
the ring of integers of discriminant $d$). For $\disc>0$ we define the
Hurwitz's class number $H(\mathfrak{d})$ by 
\begin{equation}
H(\disc) := \sum_{df^2=-\disc} \frac{h(d)}{u(d)} 
\end{equation}
if $\disc$ is a discriminant and by zero if not. A short table of the
non-zero values is given by:

$$\begin{tabular}{lcccr}
\begin{tabular}{c|c}
$\disc$ & H($\disc$)\\
\hline
3 & 1/3\\
4 & 1/2\\
7 & 1 \\
8 & 1
\end{tabular}
& & & &
\begin{tabular}{c|c}
$\disc$ & H($\disc$)\\
\hline
11 & 1\\
12 & 4/3\\
15 & 2\\
16 & 3/2
\end{tabular}
\end{tabular}$$

If $-\disc$ is a discriminant we denote $\scrO_{-\disc}$ the order of
discriminant $\disc$ in the imaginary quadratic field
$\Q[\sqrt{-\disc}]$. For $p \in \Z$ prime we define $H_p(\disc)$
to be the modified invariant as follows: 
\begin{equation}
H_p(\disc) = \left \{ \begin{array}{cl}
  0 & $ if $-\disc$ is not a discriminant$\\
  0 & $ if $p$ splits in $\scrO_{-\disc}\\
H(\disc) & $ if $p$ is inert in $\scrO_{-\disc}\\
\frac{1}{2}H(\disc) & $ if $p$ is ramified in $\scrO_{-\disc}$ but does not
divide$\\
& $ the conductor of $\scrO_{-\disc}\\
H_p(\disc/p^2)& $ if $p$ divides the conductor of $\scrO_{-\disc} \end{array}
\right.
\end{equation}
The number of embeddings of $\scrO_{-\disc}$ into any $R_i$ ($i=1, \ldots , n$)
modulo conjugation by $R_i^{\times}/\{ \pm1 \}$ is  $H_{|N|}(\disc)$ (see
\cite{Gross2} the proof of Proposition 1.9, page 122).

\noindent We want to compute the number of embeddings of
$\Z[\sqrt{N}]$ into any $R_i$, i.e. choose $\disc=4|N|$, then :
\begin{equation}
H_{|N|}(4N) = \left \{ \begin{array}{cl}
\frac{1}{2}h(4N) & \text{if } N \equiv 1 \bmod 4\\
h(N)& \text{if } N \equiv 7 \bmod 8\\
2h(N) & \text{if } N \equiv 3 \bmod 8 $ and $ N \geq 11 \end{array} \right.
\end{equation}

\noindent Note that in the case $\disc=4|N|$ an order $R_i$ on $\Ta$ appears
twice as a right order if and only if it has no embedding of
$\scrO_{4N}$. In this case it does not contribute to the sum, and
hence the number of embeddings of $\Z[\sqrt{N}]$ into the $t$ orders
in $\Ta$ is also $H_N(4N)$. With this we proved:
\begin{prop}
The number of non-equivalent Siegel points constructed is at most
$H_N(4N) t$.
\end{prop}

\begin{prop} Let $B$ be a quaternion algebra over a commutative field
$K$, and let $B_0 := \{ \beta \in B \, | \, \Tr(\beta)=0 \}$. If $\psi
: B_0 \rightarrow B_0$ is an isometry of $K$-vector spaces then there
exists an element $\beta \in B^{\star}$ such that $\sigma(x) = \beta x
\beta^{-1}$ or $\sigma(x) = - \beta x \beta^{-1} = \beta \bar x
\beta^{-1}$.
\label{Cartan}
\end{prop}

\noindent {\bf Proof.} See \cite{Vigneras} Theorem 3.3, page 12 $\square$

\begin{lem} Let $\psi : B \rightarrow B$ be an isomorphism of
$\Q$-vector spaces (respectively $\sigma: B_q \rightarrow B_q$ an
isomorphism of $\Q_q$-vector spaces) such that $\sigma(1) = 1$ and
$\sigma$ is an isometry. Then there exists an $\alpha \in B^{\star}$
(respectively $\alpha \in B_q^{\star}$) such that $\sigma(x) = \alpha x
\alpha^{-1}$ or $\sigma(x) = \alpha \bar x \alpha^{-1}$.
\label{Skolem}
\end{lem}

\noindent {\bf Proof.} Since $\sigma(1)=1$ and $\sigma$ is a morphism,
$\sigma(\Q) = \Q$. Denoting $B_0$ the trace zero elements, 
$\sigma(B_0) = B_0$ and $\sigma|_{B_0}:B_0 \rightarrow B_0$ is an
isometry. By proposition \ref{Cartan} we get two different cases:

\begin{enumerate}
\item $\sigma_{B_0}(x) = \alpha \bar x \alpha^{-1}$ for some $\alpha
\in B^{\star}$. Then $\sigma$ is the antiautomorphism given by
$\sigma(x) = \alpha \bar x \alpha^{-1}$.
\item $\sigma_{B_0}(x) = \alpha x \alpha^{-1}$ for some $\alpha \in
B^{\star}$. Then $\sigma$ is an automorphism given by $\sigma(x) =
\alpha x \alpha^{-1}$. $\square$
\end{enumerate}

\begin{thm} 
The $H_N(4N)t$ Siegel points $\{(\phi,R_i)\}$ constructed above are
\break non-equivalent.
\end{thm}

\noindent{\bf Proof.} The proof breaks in two steps. First
we will prove that for a fixed embedding of $\Z[\sqrt{N}]$ into $R$
(say $u$ is the image of $\sqrt{N}$), the $t$ left $R$-ideals give
non-equivalent points $(P,J,U)$ where $U$ is multiplication by
$u/\sqrt{|N|}$. Then we will prove that different embeddings give
non-equivalent Siegel points.

Let $I_1,I_2$ two left $R$-ideals. Abusing notation we will denote
$P_i$ the symmetric form $P_{I_i}$ and analogously for $J_i$. Suppose
there exists $W:V \rightarrow V$ an isomorphism making the Siegel
points $(P_1,J_1,U)$ and $(P_2,J_2,U)$ equivalent. Let $\beta =
W(1)$, $\sigma$ the map $\sigma(v) = W(v) \beta^{-1}$
and $V_0$ the space of elements in $V$ with trace zero. We claim that
$\sigma$ is an isometry.\\
By hypothesis $W^*P_1 = P_2$ then evaluating at $(1,1)$ we have
$$(W^*P_1)(1,1) = P_2(1,1) = \frac{2}{\Norm(I_2)}$$
By definition, $(W^*P_1)(1,1)=\frac{\Tr(W(1),\overline{W(1)})}{\Norm(I_1)} =
2\frac{\Norm(\beta)}{\Norm(I_1)}$ hence $\Norm(\beta)= \frac{\Norm(I_1)}{\Norm(I_2)}$.
Then $\|x\|/\sqrt{N} =  P_2(x,x)\Norm(I_2)/2 =
W^*(P_1(x,x))\Norm(I_2)/2  = \frac{\|W(x)\|}{\Norm(I_1)}\Norm(I_2) =
\frac{\|W(x)\|}{\|\beta\|}= \|\sigma(x)\|/\sqrt{N}$, i.e. $\sigma$ is an
isometry. Since $\sigma$ is an isometry and $\sigma(1)=1$, by lemma \ref{Skolem}
we have two different cases:

\begin{enumerate}
\item $\sigma(x) = \alpha \bar x \alpha^{-1}$ for some $\alpha \in B^\times$, i.e. $\sigma$ is an
antiautomorphism and $W(x) = \alpha \bar x \alpha^{-1} \beta^{-1}$.
\item $\sigma(x) = \alpha x \alpha^{-1}$ for some $\alpha \in B^\times$ and
$W(x) = \alpha x \alpha^{-1} \beta^{-1}$.
\end{enumerate}
We know that $W$ preserves the complex multiplication, i.e. $W^{-1}
\circ U \circ W(x) = U(x)$.\\ 
In the first case, $W^{-1}(x) = \alpha^{-1} \bar{\beta} \bar
x \alpha$. Then $W^*U(x) = W^{-1}(u \alpha \bar x \alpha^{-1}
\beta^{-1}) = \alpha^{-1} \bar \beta \bar \beta^{-1} \bar \alpha^{-1}
x \bar \alpha \bar u \alpha = x \alpha^{-1} \bar u \alpha$.  It must be
the case that $ux = x \alpha^{-1} \bar u \alpha$ for all $x \in
B$ (which is the same as saying that $ux\alpha^{-1} = x\alpha^{-1}
\bar u$) which would imply that $u \in \Q$ and is not the case. Then
we must be in the second case. 

Since $W(I_1)=I_2$, $I_2 = \alpha I_1 \alpha^{-1} \beta^{-1}$. In
particular $\alpha R \alpha^{-1} = R$, i.e. $\alpha \in Norm(R)$. Then
$I_1$ and $I_2$ have the same right order and represent the same class
between the $t$ left $R$-ideals we started with.

Assume that there is a left $R$-ideal $I$ and a left $R'$-ideal $I'$
such that $R$ and $R'$ are non-conjugate maximal orders and the Siegel
points $(P,J,U)_I$ and $(P',J',U')_{I'}$ are equivalent.  Then there exist an
isomorphism $W:V \rightarrow V$ that sends one point to the
other. Arguing as before we get the same two possible
cases for $W$. In the first case, since $W^*U = U'$ we would get that
$u'x\alpha^{-1} = x\alpha^{-1} \bar u$ for all $x \in V$. Taking $x=
\alpha$ we would get that $u' = \bar u$ and it commutes with all
elements of $V$, then it is rational which is not the case.\\ 
Then $W(x) = \alpha x \alpha^{-1} \beta^{-1}$ and $I' = \alpha I
\alpha^{-1} \beta$. In particular the orders $R$ and $R'$ are
conjugate which is a contradiction. $\square$

\subsection{Ideals associated to Siegel points}

For finding relations between the numbers $n_{[\A],[\B],\bar{\D}}$, we
 will assign to each point $z_{\A \bar \D} Q_\B$ on the Siegel space
 $\Siegel_2$ a rank $4$ $\Z$-lattice $I_z \in B$ and a basis of it such
 that the Siegel point  $(P,J,U)_I$ on this basis is $z_{\A \bar \D}
 Q_\B$. We will then prove that the left order of $I_z$ is a maximal
 order $O_{[\A],[\D]}$ with an embedding of $\Z[\sqrt{N}]$ into
 it. This will imply that the number of different values (up to a sign) for
 $n_{[\A],[\B],\bar \D}$ is at most $h(\scrO_N)^2t$ .

\begin{prop} There exists $u$ and $v$ in $B$ such
that:
\label{embeds}
\begin{itemize}
\item $\Tr(u \bar v) = 0$, $\Tr(u)=0$ and $\Tr(v)=0$
\item $\Norm(u) = |N|$
\item $\Norm(v) = |D|$
\item $u$ and $v$ are in a maximal order $R$ of $B$
\end{itemize}
\end{prop}

\noindent {\bf Proof.} Since $|N| \equiv 3 \bmod 4$, we can assume $B=
(-1,N)$. Choosing $u = j$ it is clear that $\Tr(u)=0$ and $\Norm(u)=|N|$,
hence we are looking for $v$ in $B$ such that $\Tr(uv) = 0$, $\Tr(v)
= 0$ and $\Norm(v)=|D|$. This conditions forces $v$ to have the form $v = xi +
yk$ and we are looking for an integer solution of the quadratic equation:
\begin{equation}
x^2 + |N|y^2 - |D| z^2 = 0
\label{embed}
\end{equation} 
We can assume that the solution is primitive
(i.e. $gcd(x,y,z)=1$). If $(x,y,z)$ is a solution, clearly
$gcd(z,N)=1=gcd(x,N)$ and $gcd(x,D)=1=gcd(y,D)$.\\
To prove the existence of such a solution we use the Hasse-Minkowski
principle. Clearly (\ref{embed}) has a non-zero solution over $\R$, so we need to
prove the existence of local non-zero solutions for all primes.
We consider the different cases:
\begin{itemize}
\item For a prime $p \neq N$ and $p \neq D$ the quadratic form clearly
has a local solution (see \cite{Serre} corollary 2, page 6).
\item For the prime $|N|$ by Hensel's Lemma it is enough to
look for solutions of (\ref{embed}) modulo $|N|$:
$$ x^2 - |D| z^2 \equiv 0 \bmod {|N|} \text{ iff } {(x z^{-1})}^2 \equiv |D|
\bmod{|N|}$$
This equation has solution if and only if $\left( \frac{|D|}{|N|}
\right) =1$. By the quadratic reciprocity law and the fact that $|N|
\equiv 3 \bmod 4$ this last condition is equivalent to asking that $|D|$
splits in $\Q(\sqrt{N})$ which is the case.
\item For the prime $|D|$, looking at (\ref{embed}) modulo $|D|$:
$$x^2 + |N| y^2 \equiv 0 \bmod |D| \text{ iff } N \equiv {(xy^{-1})}^2 \bmod
|D| \text{ iff } \left( \frac{N}{|D|} \right) =1$$
Which is the case since $|D|$ splits in $\Q(\sqrt{N})$.
\end{itemize}
Given $u$ and $v$ as before, consider the rank $4$ $\Z$-lattice $R =
\langle 1 , u, v, uv \rangle$. It is easy to see that $R$ is actually
an order, hence contained in a maximal one. $\square$

\noindent {\bf Remark:} if we define $R = \langle 1, \frac{1+j}{2}, v,
\left(\frac{1+j}{2}\right) v \rangle$ it is easy to see that this is
also an order. The advantage of this order is that it contains an
embedding of the ring of integers of $\Q(\sqrt{N})$, but is not maximal.\\
Let $z_{\A \bar \D} Q_{\B} = (\frac{b_1 +\sqrt{N}}{2 a_1|D|} ) \left(
\begin{array}{cc} 2a&b\\ b&2c \end{array} \right)$ , $u$ and $v$
as in proposition \ref{embeds} (choosing $u=j$). Define
\begin{equation}
I_z := \left \langle \left(\frac{b_1 - j}{2 a_1|D|}\right)av , \left(\frac{b_1
- j}{2 a_1|D|}\right)\left(\frac{|D|+bv}{2}\right), \frac{v-b}{2}, a \right \rangle
\label{ideal-def}
\end{equation}
If we denote $\phi$ the embedding of $\Q(\sqrt{N})$ into $B$ and
$\psi$ the embedding of $\Q(\sqrt{D})$ into $B$ with $\phi(\sqrt{N})=u$
and $\psi(\sqrt{D}) = v$ and choosing the basis $\B = \langle v_1,v_2
\rangle$ (where in our notation $v_1 = \frac{b+\sqrt{D}}{2}$ and $v_2
= a$) then the ideal $I_z$ was defined by:
$$I_z=\left \langle \phi\left(\frac{b_1-\sqrt{N}}{2a_1|D|} \right)
\psi(\sqrt{D}) \psi(\bar{v_2}), \phi\left(\frac{b_1-\sqrt{N}}{2a_1|D|}
\right) \psi(\sqrt{D}) \psi(\bar{v_1}), \psi (\bar{v_1}),
\psi(\bar{v_2}) \right \rangle$$
If we forget the specific basis, and think of $I_z$ just as a rank $4$
$\Z$-lattice in $B$ it is given by $I_z = \left \langle \phi
\left(\frac{b_1-\sqrt{N}}{2a_1|D|} \right) \psi(\sqrt{D}) \psi(\bar
\B) , \psi(\bar \B) \right \rangle$.
\begin{prop}
The element $\frac{1+j}{2}$ is in the left order of $I_z$.
\label{gooddefi}
\end{prop}

\noindent{\bf Proof.} This is an easy but tedious computation. We will
just give the coordinates of the product of $\frac{1+j}{2}$ with each
element of the basis of $I_z$ (given above) as a linear combination.

\noindent $\bullet$ $\left(\frac{1+j}{2}\right) a = [ba_1,-2aa_1,0,\frac{b_1+1}{2}]$.\\
\noindent$\bullet$ $\left(\frac{1+j}{2}\right) \left(\frac{v-b}{2}\right) =
[-2ca_1,ba_1,\frac{b_1+1}{2},0]$.\\
\noindent $\bullet$ $\left(\frac{1+j}{2}\right) \left(\frac{b_1 - j}{2
a_1|D|}\right)av = [\frac{1-b_1}{2},0,2ac_1,bc_1]$.\\
\noindent $\bullet$ $\left(\frac{1+j}{2}\right)\left(\frac{b_1
- j}{2 a_1|D|}\right)\left(\frac{|D|+bv}{2}\right)= [0,\frac{1-b_1}{2},bc_1,2cc_1]$. $\square$

\begin{prop} The element $a_1v$ is in the left order of $I_z$.
\end{prop}
\noindent {\bf Proof.} Since $\B$ is an ideal, it is clear that $v
\langle w_3,w_4 \rangle \subset \langle w_3,w_4 \rangle$. By the way
we choose $v$, it satisfies $vj=-jv$, then 
\begin{equation}
(a_1v)\left(\frac{b_1-j}{2a_1|D|} \right) =
\left(\frac{b_1-j}{2a_1|D|}\right) (-a_1v) +\frac{b_1}{|D|} v
\label{lor}
\end{equation}
For the part corresponding to the first two elements of $I_z$ note
that they can be written as $\left(\frac{b_1-j}{2a_1|D|}\right) v (a)$ and
$\left(\frac{b_1-j}{2a_1|D|}\right) v
\left(\frac{v-b}{2}\right)$. Since $\B$ is an ideal, $v \B \subset \B$
and the assertion follows from equation (\ref{lor}). $\square$ 

\begin{cor} The order $R = \langle 1, \frac{1+j}{2},a_1v,\frac{1+j}{2}a_1v \rangle$
is contained in the left order of $I_z$ and has discriminant
$(a_1^2ND)^2$ or index $a_1^2|D|$ in a maximal order. 
\label{Rorder}
\end{cor}

\noindent {\bf Proof.} It is clear that $R$ is in the left order of
$I_z$ by the previous two propositions. It is also clear that it is an
order. To compute its discriminant, note that the bilinear matrix
associated to it is:
\begin{center}
$\left( \begin{array}{cccc}
2&1&0&0\\
1&\frac{1-N}{2}&0&0\\
0&0&2a_1^2|D|&a_1^2|D|\\
0&0&a_1^2|D|&a_1^2|D|\frac{1-N}{2}
\end{array} \right)$\\
\end{center}
Then note that the index in a maximal order (which has discriminant
$N^2$) is the square root of the discriminant. $\square$

\begin{thm} Let $U$ be the complex multiplication associated to
$\frac{-j}{\sqrt{|N|}}$. Then the Siegel point $(P,J,U)_{I_z}$ associated to
the ideal $I_z$ in the given basis is $z_{\A \bar \D} Q_{\B}$. 
\label{ideal}
\end{thm}

\noindent {\bf Proof.} This is a straightforward computation so we
omit the details. Just check that the given basis of $I_z$ is symplectic,
i.e. that the matrix $J(x,y)$ in the given basis is a multiple of the
matrix $\left(
\begin{array}{cc} 0& I_2\\ -I_2&0 \end{array} \right)$ (since $J(x,y)$
is skewsymmetric there are half the conditions to check), and that
the matrix $U$ associated to the point $z_{\A \bar \D} Q_{\B}$ is the
same as the complex multiplication matrix on $I_z$. $\square$

\begin{thm}
The lattice $I_z$ is an ideal for a maximal order.
\end{thm}

\noindent {\bf Proof.} The strategy is to prove that the quadratic form associated to
the ideal $I_z$ is locally equivalent to the maximal order one for all
primes. We need the next lemma:

\begin{lem} The quadratic form associated to the lattice
$I_z$ has discriminant $N^2$.
\end{lem}

\noindent {\bf Proof.} The bilinear form is the same as the Siegel
point $z_{\A \bar\D} Q_\B$ hence its bilinear form matrix is:
\begin{center}
$B_I=\left( \begin{array}{cc}
2c_1 Q_\B& b_1I_2\\
b_1I_2& 2a_1DQ_\B^{-1} \end{array} \right)$\\
\end{center}
Since $Q_\B$ has determinant $D$, it is an easy computation
to prove that the determinant of this matrix is $N^2$ (using that
$b_1^2-4a_1c_1|D|=N$). $\square$

A maximal order for $B=(-1,N)$ is given by $O=\langle \frac{1+j}{2},\frac{i+k}{2},j,k
\rangle$ (see Proposition 5.2, page 369 of \cite{Pizer}), then it is
easy to compute the matrix of the quadratic form trace and to check
that it has discriminant $N^2$, and is an improperly
primitive integral form. Since the discriminant of both forms is a
unit for all primes $p \neq |N|$ then they are locally equivalent (see
Corollary of Theorem 3.1 of \cite{Cassels}, page 116).
Hence $(I_z)_p$ is locally principal for all primes $p \neq |N|$.\\
For the ramified prime, $D(I_z)=N^2$ hence it is locally
principal. Locally principal ideals have the same discriminant as
their left orders hence $O_l(I_z)$ is maximal. $\square$

\subsection{Comparing Siegel Points}

If $I$ is an ideal for a maximal order, and $U$ a complex
multiplication, the Siegel point associated to $(U,I)$ is the same as
the one associated to the point $(U,I\alpha)$ for any $\alpha \in
B^\times$ (with the same choice of basis). Suppose two Siegel points
$z$ and $z'$ have equivalent ideals $I_z$ and $I_{z'}$, say $I_z =
I_{z'} \alpha$ for some $\alpha \in B^\times$. Then since the complex
multiplication is the same for all the ideals we constructed, the two
Siegel points are equivalent by proposition \ref{equivalence}. Let $M$
be the matrix in $Sp_4(\Z)$ making the change of basis between 
$I_z$ and $I_{z'}\alpha$.

\begin{lem} .The matrix $M$ is in the subgroup $\Gamma_{1,2}$.
\label{equivideals}
\end{lem}

\noindent {\bf Proof.} Let $M=\left(
\begin{array}{cc}
A&B\\
C&D \end{array} \right)$,  $z
= \left( \frac{b_1+\sqrt{N}}{2a_1} \right) Q$ and $z' = \left(
\frac{b_1'+\sqrt{N}}{2a_1'} \right) Q'$ where $Q$ and $Q'$ have even
diagonal. Since $M$ sends the bilinear
form associated to the ideal $I_z$ to the bilinear form associated to
the ideal $I_{z'}\alpha$, 
$$\left(\begin{array}{cc}
A & B\\
C & D \end{array} \right)^t 
\left( \begin{array}{cc}
2c_2Q & b_2I_2\\
b_2 I_2& 2 a_2 Q^{-1} \end{array} \right)
\left(\begin{array}{cc}
A & B\\
C & D \end{array} \right) = 
\left( \begin{array}{cc}
2c_2' Q & b_2'I_2\\
b_2' I_2& 2 a_2' Q'^{-1} \end{array} \right)$$
By the way we choose generators, $b_i \equiv 1 \bmod 4$, $i=1,2$ (also $b_i'
\equiv 1 \bmod 4$, $i=1,2$) hence $2Q \equiv 
\left( \begin{array}{cc}
0&2\\
2&0 \end{array} \right) \bmod 4$. Let $J := 
\left( \begin{array}{cc}
0&1\\
1&0 \end{array} \right)$. Looking at the first $2 \times 2$
matrix of the previous equality modulo $4$ we get:
$2c_2 A^t J A + C^tA+A^tC + 2a_2 C^tJC \equiv 2J \bmod 4$. In
particular $4$ divides the diagonal. \\
If $A := \left(
\begin{array}{cc}
a&b\\
c&d \end{array} \right)$ then $A^tJA = \left(
\begin{array}{cc}
2ac&ad+bc\\ ad+bc&2bd \end{array} \right)$ hence $4$ divides the
diagonal of $2c_2A^tJA$ and $2a_2C^tJC$. Also $A^tC$ is symmetric
hence $A^tC+C^tA = 2A^tC$ and we get that $2$ divides the diagonal of
$A^tC$. The proof for $B^tD$ is analogous looking at the last $2
\times 2$ matrix. $\square$

\begin{prop} For fixed ideals $\A$ and $\D$, the left order of
$I_{z_{\A \D} Q_{\B}}$ is independent of the ideal $\B$.
\label{changingB}
\end{prop}

\noindent {\bf Proof.} We know $I_z = \left \langle \phi
\left(\frac{b_1-\sqrt{N}}{2a_1|D|} \right) \psi(\sqrt{D}) \psi(\bar
\B) , \psi(\bar \B) \right \rangle$. The ideal $\B_q :=\B \otimes
\Z_q$ is principal, hence there exists an element $\delta_q \in
L_q:=\Q_q(\sqrt{D})$ such that $\B_q = \scrO_L \delta_q$. Then 
$I_z \otimes \Z_q = \left \langle \phi\left(\frac{b_1+\sqrt{N}}{2a_1|D|}
\right) \psi(\sqrt{D}) \psi(\scrO_L), \psi(\scrO_L) \right \rangle
\bar \delta_q$, hence its left order is clearly independent of $\B$. $\square$

\begin{prop} Let $\A$ and $\A'$ two equivalent ideals of
  $\scrO_K$ prime to $\D$, say $\A' = \alpha \A$. Then $\phi(\alpha^{-1}) I_{z_{\A \D}
  Q_{\B}} = I_{z_{\A' \D} Q_{\B}}$.
\label{sametypeorders}
\end{prop} 

\noindent {\bf Proof.} It is enough to prove that $I_{z_{\alpha \A \D}
  Q_{\B}} \subseteq \phi(\alpha^{-1}) I_{z_{\A \D} Q_{\B}}$. Then $I_{z_{\A \D}
  Q_{\B}} \subseteq \phi(\alpha) I_{z_{\alpha \A \D} Q_{\B}}$ and the result
  follows.\\
Without loss of generality we may assume that
$\A$ and $\A'$ are prime to each other, then we can choose basis such
that $\bar \A \bar \D = \langle a|D|, \frac{b-\sqrt{N}}{2} \rangle$
and $\bar \A' \bar \D = \langle a'|D|, \frac{b-\sqrt{N}}{2}
\rangle$. Then there exists $M = \left( \begin{array}{cc}
x_1&x_2\\
x_3&x_4 \end{array} \right)  \in Sl_2(\Z)$ such that
\begin{eqnarray}
\label{eqn1}
x_1 a|D| \bar \alpha + x_2 \left(\frac{b-\sqrt{N}}{2} \right) \bar
\alpha =& a' |D|\\
\label{eqn2}
x_3 a|D| \bar \alpha + x_4 \left(\frac{b-\sqrt{N}}{2} \right) \bar
\alpha =& \frac{b-\sqrt{N}}{2} 
\end{eqnarray}
{\bf claim:} $D \mid x_2$. If $\bar \alpha = \frac{1}{a}
\left(\frac{\alpha_1 + \alpha_2 \sqrt{N}}{2} \right)$ with $\alpha_i
\in \Z$ looking at the imaginary parts of the above equalities we get
that 

$\left \{
\begin{array}{cc}
\frac{\sqrt{N}}{4}\left(2 x_1 a|D| \alpha_2 + x_2 (b \alpha_2 -
\alpha_1)\right) =& 0\\ 
\frac{\sqrt{N}}{4} \left(2x_3 a|D| \alpha_2 + x_4 (b \alpha_2 -
\alpha_1)\right) =& \frac{-a\sqrt{N}}{2} 
\end{array}\right.$

this implies the claim. If $\B = \langle w_1,w_2
\rangle$, 
$\phi(\alpha^{-1}) I_{z_{\A \D}}$ = $\left \langle \phi
\left(\bar \alpha (\frac{b-\sqrt{N}}{2a|D|}) \right) \psi(\sqrt{D} \bar
w_1) , \right.$ $\left . \phi
\left(\bar \alpha (\frac{b-\sqrt{N}}{2a|D|}) \right)
\psi(\sqrt{D} \bar w_2), \phi(\alpha^{-1}) \psi(\bar w_1), \phi(\alpha^{-1})
\psi(\bar w_2)  \right \rangle$. Since $a \bar \alpha = \alpha^{-1}a'$
equation (\ref{eqn1}) implies that  $$x_3
\phi(\alpha^{-1})+ x_4
\phi\left(\alpha^{-1}(\frac{b-\sqrt{N}}{2a|D|})
\right)= \phi\left(
\frac{b-\sqrt{N}}{2a'|D|}\right) $$

Since $\B$ is an ideal, $\sqrt{D} w_i \in \B$ hence
$\phi(\frac{b-\sqrt{N}}{2a'|D|}) \psi(\sqrt{D} \bar w_i) \in
\phi(\alpha^{-1}) I_{z_{\A \D} Q_{\B}}$ for $i=1,2$. Since $D \mid x_2$ equation
(\ref{eqn2}) can be written as
$$x_1\phi(\alpha^{-1}) + \frac{x_2}{|D|} \phi\left(\alpha^{-1}
(\frac{b-\sqrt{N}}{2a|D|})\right) \psi(\sqrt{D})^2= 1$$
which implies that $\psi(\bar w_i) \in \phi(\alpha^{-1}) I_{z_{\A \D}
  Q_{\B}}$ for $i=1,2$. $\square$

\begin{cor} If $\A$, $\A'$ are two equivalent ideals in $\scrO_K$
  prime to $\D$ then the ideal $I_{z_{\alpha \A \D}
  Q_{\B}}$ and $I_{z_{\alpha \A' \D} Q_{\B}}$ have equivalent left
  orders.
\label{sametype}
\end{cor}

\begin{prop} Let $\D$ and $\D'$ be two split prime ideals of $\Q[\sqrt{N}]$
 of norms $|D|$ and $|D'|$ respectively such that $\D'=\mu \D$.
Let $\B$ and $\B'$ be ideals of $\Q[\sqrt{D}]$ and of
$\Q[\sqrt{D'}]$ respectively. Then the ideals $I_{z_{\A \D} Q_\B}$ and
$I_{z_{\A \D'}Q_{\B'}}$ have the same left order if following the notation of
proposition \ref{embeds} we take $v'= \mu v$.
\label{ppal}
\end{prop}

\noindent{\bf Proof.} We are abusing notation while stating this
theorem, since $\mu$ is an element of $\Q[\sqrt{N}]$. We will not
distinguish between an element in $B$ or in $\Q[\sqrt{N}]$ via the
identification $\sqrt{N} \mapsto j$, and the case will be clear from
the context.\\ 
By proposition \ref{changingB} it is enough to restrict to the case
$\B$ and $\B'$ principal. In this case we will prove that the ideals
associated to them are slightly different and use this to prove the
proposition. We can choose basis such that $\D = \langle |D|,
\frac{b_1+\sqrt{N}}{2} \rangle$ and $\D' = \langle |D'|, \frac{b_1
+\sqrt{N}}{2} \rangle$. Let $\mu = \frac{\alpha}{|D|}+\frac{\beta}{|D|}
\sqrt{N}$. Since $\mu \left(\frac{b_1+\sqrt{N}}{2} \right) \in \D'$
and $\mu^{-1} \left(\frac{b_1+\sqrt{N}}{2} \right) \in \D$,
$\frac{\alpha + \beta b_1}{|D|} \in \Z$ and 
$\frac{\alpha -\beta b_1}{|D'|} \in \Z$. \\
Since $b = 1$ the definition of the ideals is:
\begin{itemize}
\item $I_\D := I_{z_{\A \D} Q_\B}= \left\langle \left(\frac{b_1-j}{2a_1|D|} \right) v,
\left(\frac{b_1-j}{2a_1|D|} \right) \left(\frac{v+|D|}{2}\right),
\frac{v-1}{2},1 \right \rangle$ \\
\item $I_{\D'} := I_{z_{\A \D'} Q_\B'}=\left\langle \left(\frac{b_1-j}{2a_1|D'|} \right) v',
\left(\frac{b_1-j}{2a_1|D'|} \right) \left(\frac{v'+|D'|}{2}\right),
\frac{v'-1}{2},1 \right \rangle$
\end{itemize}
where $v$ and $v'$ are the elements of norm $|D|$ and $|D'|$
respectively as in proposition \ref{embeds}. We will write the
elements of $I_{\D'}$ in the basis of $I_{\D}$, the other case follows
from symmetry.
\begin{itemize}
\item $\frac{v'-1}{2}=[-a_1\beta,0,\frac{\alpha+b_1\beta}{|D|},\frac{\alpha+b_1\beta
+D}{2|D|}]$\\
\item $\left( \frac{b_1-j}{2a_1|D'|} \right) v'= [\frac{\alpha-\beta
    b_1}{|D'|},0,4 \beta c, 2\beta c]$  which has integer coefficients\\
\item $\left( \frac{b_1-j}{2a_1|D'|} \right) \left(
\frac{v'+|D'|}{2} \right)=[\frac{\alpha - \beta
    b_1-|D'|}{2|D'|},1,2\beta c,\beta c]$
\end{itemize}
We cannot say that the two ideals are the same, since the numbers
$\alpha$ and $\beta$ may have a $2$ in the denominator, but $(I_\D)_p
= (I_{\D'})_p$ for all primes $p \neq 2$. In particular if we denote
$O_\D$ and $O_{\D'}$ the left order of $I_\D$ and $I_{\D'}$
respectively, we get that $(O_D)_p = (O_{\D'})_p$ for all $p \neq
2$. Since the denominators are at most $2$ it is easy to check that
$4O_\D +\Z \subset O_{\D'}$, and has index at most $2^8$. By corollary
\ref{Rorder}, the order $R \subset O_{\D'}$ with index $a_1^2|D|$,
which is odd. Then $4O_D + R = O_{\D'}$. Also $4O_{\D} +R = O_\D$
hence both orders are the same. $\square$

By theorem \ref{comparison} we know that the numbers
$n_{[\A],[\B],\bar \D}$ depend (up to multiplication by $\pm1$) on the
equivalence class of $\A$, the equivalence class of $\D$ and the class
of $z_{\A \D}Q_\B$ modulo $\Gamma_{12}$. If we fix the class of $\A$
and the class of $\D$ we can associate ideals to the points $z_{\A
\D}Q_\B$ as in (\ref{ideal-def}) and by Proposition \ref{comparison}
they all have the same left order. Then by Corollary \ref{sameRorder} we
get at most $t(B)$ different points in the Siegel space. This implies:
\begin{thm} The number of different $n_{[\A],[\B],\bar \D}$ up to
  multiplication by $\pm1$  in $\M$
is at most $h(\scrO_K)^2 t(B)$, where $t(B)$ is the type number for
maximal orders.
\end{thm}

Note that this number is independent of the class number of
$\scrO_L$. With all these results we return and finish the proof of Theorem
\ref{main-teo}: 

Given $\A$ and $[\D]$ as before we associate to them a maximal order
$O_{\A,[\D]}$. For any left $O_{\A,[\D]}$-ideal $I$ we want to
define the number $m_{[\A],I}([\D])$.
\begin{itemize}
\item If there
exists a pair $(\D',\B)$ where $\D' \in \scrO_K$ is a prime ideal of
norm $D'$ congruent to $3$ modulo $4$, $\D' \sim \D$ and $\B$ is an
ideal of $\Q(\sqrt{-D'})$ such that $I=I_{z_{\A  \D'} Q_\B}$, we
define $m_{\A,I}([\D]) = \xi_2 n_{\A,[\B],\bar \D'}$.\\
The number $\xi_2$ is chosen such that
$m_{\A,I}([\D])$ is a complex number in the upper half
plane union $\R_{\geq 0}$.
\item If no such pair exists we define $m_{\A,I}([\D])=0$.
\end{itemize}
\begin{prop}
This definition is ``independent'' of the equivalent class of
the ideal $\A$. 
\end{prop}

\noindent {\bf Proof.} By Corollary \ref{sametype} if two ideal $\A$, $\A'$
are equivalent (say $\A' = \alpha \A$), their left orders are
conjugate. Furthermore a bijection between 
left $O_{\A,[\D]}$-ideals and left $O_{\A',[\D]}$-ideals is given by
multiplication on the right by $\phi(\alpha^{-1})$ (by Proposition
\ref{sametypeorders}). Since the number $n_{[\A],[\B],\bar \D'}$ is
independent of the equivalent class of $\A$ this map preserves
the numbers $\{m_{\A,I}([\D])\}$. $\square$

Hence we think of the numbers $m_{\A,I}([\D])$ as defined on
equivalence classes and denote them $m_{[\A],I}([\D])$. 

The formula (\ref{main-for4}) says:
\begin{equation*}
L(\psi_\D,1) =\frac{2\pi}{w\sqrt{|D|}} \eta(\bar\D) \eta(\scrO_K) \left( \sum_{[\A] \in Cl(\scrO_K)} 
 \sum_{[\B] \in Cl(\scrO_L)}
n_{[\A],[\B],\bar \D} \right) 
\end{equation*}
To the Siegel point $z_{\A \bar \D} Q_\B$ we associate the left
$O_{[\A],[\D]}$-ideal $I_{\B}$ as in (\ref{ideal-def}). Given $I$ a
left $O_{[\A],[\D]}$-ideal, we define
$$r(\D,[\A],I)= \left\{ \begin{array}{ll} \sum_{\{\B \in \scrO_L |
I_\B \sim I\}} n_{[\A],[\B],\bar D}/m_{\A,I}([\D]) & \text{\qquad if
$m_{\A,I}([\D]) \neq 0$}\\ 0 & \text{\qquad
otherwise}\end{array}\right.$$ 
Lemma \ref{equivideals} and Theorem \ref{comparison} imply that if the ideals
$I_\B$ and $I_{\B'}$ are equivalent, $n_{[\A],[\B],\bar D} = \pm
n_{[\A],[\B'],\bar D}$ hence $r(\D,[\A],I) \in \Z$. Rearranging the sum
we get:
\begin{equation*}
L(\psi_\D,1) =\frac{2\pi}{w\sqrt{|D|}} \eta(\bar\D) \eta(\scrO_K) \left( \sum_{[\A] \in Cl(\scrO_K)} 
 \sum_{I} r(\D,[\A],I) m_{[\A],I}([\D]) \right) 
\end{equation*}
as claimed $\square$

\noindent{\bf Question:} is it true that for any left
$O_{[\A],[\D]}$-ideal $I$ there exists a pair $(\D',\B)$ such that $I
\sim I_{z_{\A  \D'} Q_\B}$?\\
All the examples we computed show this is the case.

\begin{prop}
Let $\A$ be an ideal of $\Q(\sqrt{N})$, then $n_{[\A],[\B],\bar \D}$
and $n_{[\scrO_K],[\B],\bar \D}$ differ by a unit in a quadratic
extension of $\M$.
\end{prop}

\noindent{\bf Proof.} Let $\sigma_{\A}$ be the automorphism of $H$
corresponding to the ideal $\A$ via the Artin-Frobenius map. Then we proved that
$\left( \frac{\theta(z_{\scrO_K \D}Q_\B)}{\eta(\D)\eta(\scrO_K)}
\right) ^{\sigma_A} = \frac{\theta(z_{\A \D}Q_\B)}{\eta(\A \D)
\eta(\A)}$. Hence $n_{[\A],[\B],\bar \D} = \left(\frac{\eta(\A)
\eta(\A \D)}{\eta(\D)\eta(\scrO_K)\psi_{\bar \D}(\A)}\right)
(n_{[\scrO_K],[\B],\bar \D})^{\sigma_\A}$. Note that the quotient of
etas squared is in $H$ while $\psi_{\bar \D}(\A)$ is in $T$, hence
$\zeta:=\left(\frac{\eta(\A) \eta(\A
\D)}{\eta(\D)\eta(\scrO_K)\psi_{\bar \D}(\A)}\right)$ is in a
quadratic extension of $\M$. Clearly $\Norm(\zeta) = 1$ as
required. $\square$

\section{The class number one case}

We study now the case of imaginary quadratic fields with
class number equal to one. In this case $n_{[\A],[\B],\bar \D}$ are
rational integers for any choice of $\D$.  There are just six such
cases (we exclude the case $N=-3$) so we can study all this cases by
numerical computations. Here are some examples:
\subsection{Case $N=-7$}

This case is the easiest one since the class number in the
quaternion algebra is also one. Then the numbers $n_{[\A],[\B],\bar
\D}$ are integers and differ by a unit. 

\begin{thm}
Let $N=-7$ and $\D$ be any ideal of prime norm congruent to $3$ modulo
$4$. Then $L(\psi_\D,1) \not = 0$.
\label{non-vanishing}
\end{thm}
\noindent{\bf Proof.} By proposition \ref{conjugationofB} we know that
the number associated to an ideal $\B$ is the same as the one
associated to $\bar \B$. For a prime ideal $\D$ let $\Omega =
\eta(\bar \D) \eta (\scrO_K) \frac{2\pi}{w\sqrt{|D|}}$ where
$-D=\Norm(\D)$ and $w$ is the number of units in $\Q[\sqrt{D}]$. 
The formula \ref{main-for4} for $L(\psi,1)$ reads:
\begin{equation}
L(\psi,1)= \left(\sum_{[\B] \in Cl(\scrO_L)} n_{[\scrO_K],[\B],\bar \D}\right)
\Omega = \left(n_{[\scrO_K],[\scrO_L],\bar \D} + 2 {\sum_{[\B]\in \Phi}} n_{[\scrO_K],[\B],\bar \D}\right)\Omega
\end{equation}
where $\Phi$ is a maximal subset of $Cl(\scrO_L)$ such that $[\scrO_L]
 \not \in \Phi$ and if $[\B] \in \Phi$ then $[\bar \B] \not \in \Phi$.

Taking the maximal order $O$ as left $O$-ideal representative, we see
that the number associated to it is $1$ up to a sign, then
$\frac{L(\psi,1)}{\Omega} \equiv 1 \bmod 2$. $\square$

In the next table, we list some of the numbers $n_{[\scrO_K],[\B],\bar
\D}$ to show the behavior of the sign.

\begin{center} 
\begin{tabular}{|c|c|c|} \hline
\bf $D$ & $\B$ & $n_{[\A],[\B],\bar\D}$ \\
\hline
11 & [1,-1,3] & 1\\
23 & [1, -1, 6] & 1\\
23 & [13, -17, 6] & -1\\
23 & [13, 17, 6] & -1\\
43 & [1, -1, 11] & -1\\
67 & [1, -1, 17] & 1\\

\hline
\end{tabular}

\begin{tabular}{|c|c|c|c|} 
\hline
\bf $D$ & $\B$ & $n_{[\A],[\B],\bar\D}$\\
\hline

71 & [1, -1, 18] & -1\\
71 & [19, 9, 2] & -1\\
71 & [19, -9, 2] & -1\\
71 & [29, 33, 10] & 1\\
71 & [29, -33, 10] & 1\\
71 & [43, 141, 116] & -1\\
71 & [43, -141, 116] & -1\\
\hline
\end{tabular}
\end{center}

\subsection{Case $N=-11$}

In this case the quaternion algebra has type number $2$ for maximal
orders, so we get two
different integers associated to different $\D$'s. Each number
$n_{[\scrO_K],[\B],\bar \D}$ will be
associated to an ideal class. Let $B=(-1,-11)$ be the quaternion
algebra ramified at $11$ and infinity. Let $O := \langle \frac{1}{2}+\frac{j}{2}, \frac{i}{2}+\frac{k}{2},j,k
\rangle$ be a maximal order and $I$ a non-principal ideal. Here is a
table of $n_{[\scrO_K],[\B],\bar \D}$ for different values of
$D$ and $\B$, writing down the associated ideal also.

\begin{center}
\begin{tabular}{|c|c|c|c|} \hline
\bf $D$ & $\B$ & $n_{[\A],[\B],\bar\D}$ & Ideal \\
\hline

23 & [1, -1, 6] & 2 & $I_1$\\
23 & [13, -17, 6] & 0 & O\\
23 & [13, 17, 6] & 0 & O\\
31 & [1, -1, 8] & -2 & $I_1$\\
31 & [5, 17, 16] & 0 & O\\
31 & [5, -17, 16] & 0 & O\\
47 & [1, -1, 12] & 0 & O\\
47 & [7, -17, 12] & 2 & $I_1$\\
47 & [7, 17, 12] & 2 & $I_1$\\
47 & [17, -53, 42] & 0 & O\\
47 & [17, 53, 42] & 0 & O\\
\hline
\end{tabular}
\end{center}
Note that the number $0$ is associated to the principal ideal, while
the number $2$ is associated to $I_1$. With the same reasoning as in
theorem \ref{non-vanishing} we can get a partial result proving that
the ideals $\D$ such that $z_{\D}Q_{\scrO_L}$ is associated to the
ideal $I_1$ have a non-vanishing L-series.\\
Following the method described in \cite{Pacetti}, taking $\{O,I_1\}$
as representatives for the maximal order and constructing the Brandt
matrices for level $11^2$ we get that the eigenvector associated to
the modular form of weight $2$ and level $11^2$ is
$[0,0,0,1,-1,0,0,0,1,-1]$. The first three zeros correspond to the
principal ideal, and the $\pm 1$ to $I_1$. Then the number associated
to each ideal is the same as the one associated to it via
$n_{[\scrO_K],[\B],\bar \D}$, since the eigenvector is well defined up
to a constant.

\subsection{Case $N=-163$}

Let $B=(-1,-163)$ be the quaternion algebra ramified at $163$ and
infinity. In this case, the class number for maximal orders is $14$
while the type number is $8$. Consider the maximal order 
$O:=\langle 1, i,
\frac{1}{2}+\frac{j}{2},\frac{i}{2}+\frac{k}{2}\rangle$
A set of representatives of left $O$-ideals is
given by $\{I_j\}_{j=1}^{14}$ with $I_1=O$ and
\begin{itemize}
\item $I_2:=\langle2,2i,\frac{1}{2}+i+\frac{j}{2},-1+\frac{i}{2}+
\frac{k}{2} \rangle$
\item $I_3:= \langle 3, 3i, \frac{1}{2}+i+\frac{j}{2},
-1+\frac{i}{2}+\frac{k}{2} \rangle$
\item $I_4:=\langle 3, 3i, \frac{-1}{2}+i+\frac{j}{2},
-1-\frac{i}{2}+\frac{k}{2} \rangle$
\item $I_5:= \langle 6, 6i, \frac{1}{2}+i+\frac{j}{2},
-1+\frac{i}{2}+\frac{k}{2} \rangle$
\item $I_6 := \langle 6, 6i, \frac{-1}{2}+i+\frac{j}{2},
-1-\frac{i}{2}+\frac{k}{2} \rangle$
\item $I_7:= \langle 4, 4i, \frac{3}{2}+i+\frac{j}{2},
-1+\frac{3i}{2}+\frac{k}{2} \rangle $
\item $I_8:= \langle 4, 4i, \frac{-3}{2}+i+\frac{j}{2},
-1-\frac{3i}{2}+\frac{k}{2} \rangle$
\item $I_9:= \langle 6, 6i, \frac{5}{2}+i+\frac{j}{2},
-1+\frac{5i}{2}+\frac{k}{2} \rangle$
\item $I_{10} := \langle 6, 6i,\frac{-5}{2}+i+\frac{j}{2},
-1-\frac{5i}{2}+\frac{k}{2} \rangle$
\item $I_{11}:= 5, 5i, \frac{1}{3}+2i+\frac{j}{2},
-2+\frac{i}{2}+\frac{k}{2} \rangle$
\item $I_{12}:= \langle 5, 5i, \frac{-1}{2}+2i+\frac{j}{2},
-2-\frac{i}{2}+\frac{k}{2} \rangle$
\item $I_{13}:= \langle 7, 7i, \frac{5}{2}+3i+\frac{j}{2},
-3+\frac{5i}{2}+\frac{k}{2} \rangle$
\item $I_{14}:= \langle 7, 7i, \frac{-5}{2}+3i+\frac{j}{2},
-3-\frac{i}{2}+\frac{k}{2} \rangle$
\end{itemize}
The pairs of ideals $(I_{2j+1},I_{2j+2})$ with $j=1,\ldots,6$, have the
same right order, hence each pair will have the same integer associated. For
the table we consider the range of primes between $150$ and $200$ so
as to get all the ideals $\{I_j\}$ associated to some number
$n_{[\scrO_K],[\B],\bar \D}$. The table is:
\begin{center}
\begin{tabular}{|c|c|c|c|} \hline
\bf $D$ & $\B$ & $n_{[\A],[\B],\bar\D}$ & Ideal \\
\hline

151 & [1, -1, 38] & 20 & $I_2$\\
151 & [29, 9, 2] & 14 & $I_8$\\
151 & [29, -9, 2] & 14 & $I_8$\\

\hline
\end{tabular}

\begin{tabular}{|c|c|c|c|} 
\hline
\bf $D$ & $\B$ & $n_{[\A],[\B],\bar\D}$ & Ideal \\
\hline
151 & [11, -5, 4] & 8 & $I_{13}$\\
151 & [11, 5, 4] & 8 & $I_{14}$\\
151 & [43, 137, 110] & 4 & $I_{12}$\\
151 & [43, -137, 110] & 4 & $I_{12}$\\
167 & [1, -1, 42] & 0 & $I_1$\\
167 & [157, 33, 2] & -20 & $I_2$\\
167 & [157, -33, 2] & -20 & $I_2$\\
167 & [61, 65, 18] & -2 & $I_4$\\
167 & [61, -65, 18] & -2 & $I_3$\\
167 & [29, 93, 76] & -10 & $I_6$\\
167 & [29, -93, 76] & -10 & $I_5$\\
167 & [127, -177, 62] & -14 & $I_7$\\
167 & [127, 177, 62] & -14 & $I_8$\\
167 & [19, -21, 8] & -12 & $I_9$\\
167 & [19, 21, 8] & -12 & $I_{10}$\\
179 & [1, -1, 45] & 0 & $I_1$\\
179 & [19, 45, 29] & 2 & $I_3$\\
179 & [19, -45, 29] & 2 & $I_4$\\
179 & [13, 17, 9] & 4 & $I_{12}$\\
179 & [13, -17, 9] & 4 & $I_{11}$\\
199 & [1, -1, 50] & 0 & $I_1$\\
199 & [31, -69, 40] & -20 & $I_2$\\
199 & [31, 69, 40] & -20 & $I_2$\\
199 & [43, -133, 104] & -4  & $I_{12}$\\
199 & [43, 133, 104] & -4  & $I_{11}$\\
199 & [13, 29, 20] & -14 & $I_8$\\
199 & [13, -29, 20] & -14 & $I_7$\\
199 & [131, 453, 392] & -8  & $I_{14}$\\
199 & [131, -453, 392] & -8  & $I_{13}$\\
\hline
\end{tabular}
\end{center}
The eigenvector for the Brandt matrices corresponding to the form of
weight $2$ and level $167^2$ is given by the vector $[0,10, 1, 1, 5, -5,
7, -7, -6, 6, 2, 2, -4, 4]$ with respect to the maximal order
representatives $\{I_j\}$.\\
Considering all the class number 1 imaginary quadratic fields (the
computations being the same in all cases), we can prove:

\begin{thm}
Let $E$ be a CM elliptic curve over $\Q$ of level $p^2$. Then the
coordinate of the eigenvector of the Brandt matrices associated to
$E$ on the place corresponding to an ideal $I$ is given up to a sign
by $m_{[\scrO_K],I}([\D])$. 
\end{thm}

\end{document}